\documentclass[preprint,review,11pt]{elsarticle}

\usepackage{natbib}
 \bibpunct[, ]{(}{)}{,}{a}{}{,}%
\usepackage{comment}
\usepackage{subfigure}
\usepackage{floatrow}
\usepackage{lscape}
\usepackage{newunicodechar}
\newunicodechar{✓}{\checkmark}
\usepackage{enumerate}
\usepackage{mathrsfs}
\usepackage[justification=centering]{caption}
\usepackage{algorithm,setspace}
\usepackage[noend]{algpseudocode}
\usepackage{epstopdf}
\usepackage{hyperref}
\usepackage{amsmath}
\usepackage{amssymb}
\usepackage{mathtools}
\usepackage{cases}

\makeatletter
\def\@author#1{\g@addto@macro\elsauthors{\normalsize%
    \def\baselinestretch{1}%
    \upshape\authorsep#1\unskip\textsuperscript{%
      \ifx\@fnmark\@empty\else\unskip\sep\@fnmark\let\sep=,\fi
      \ifx\@corref\@empty\else\unskip\sep\@corref\let\sep=,\fi
      }%
    \def\authorsep{\unskip,\space}%
    \global\let\@fnmark\@empty
    \global\let\@corref\@empty  
    \global\let\sep\@empty}%
    \@eadauthor={#1}
}
\makeatother

\newtheorem{theorem}{Theorem}

\newtheorem{defn}{Definition}

\newcommand{\citel}[1]{(\citeauthor{#1}, \citeyear{#1})}
\newcommand{\citenl}[1]{\citeauthor{#1}, \citeyear{#1}}

\newcommand{\mI}{\mathcal{I}}
\newcommand{\bb}{\boldsymbol{\beta}}

\journal{European Journal of Operational Research}

\begin{document}

\begin{frontmatter}



\title{On the Equity-Efficiency Trade-off in Food-Bank Network Operations}

\author{Mohammad Firouz\corref{cor1}\fnref{fn1}}
\cortext[cor1]{Corresponding Author; E-mail: mfirouz@uab.edu; Phone: (205)934-8830}

\author{Linda Li\fnref{fn2}}
\author{Barry Cobb\fnref{fn3}}
\author{Feibo Shao\fnref{fn2}}

\address[fn1]{Department of Management, Information Systems, and Quantitative Methods, University of Alabama at Birmingham, Birmingham, AL 35294.}
\address[fn2]{College of Business, Missouri State University, Springfield, MO 65897.}
\address[fn3]{Department of Economics and Business, Virginia Military Institute, Lexington, VA 24450.}

\begin{abstract}

\noindent In this paper, we present a novel modelling perspective to the food-bank donation allocation problem under equity and efficiency performance measures. Using a penalty factor in the objective function, our model explicitly accounts for both efficiency and equity, simultaneously. We give the tightest lower and upper bounds of the penalty factor, which can conveniently characterize closed-form optimal solutions for the perfect efficiency and perfect equity cases. Testing our model on the full spectrum of our penalty factor, using real data from Feeding America, we demonstrate that the solutions from our model dominate those of a benchmark from the literature in terms of equity, efficiency, and utilization equity (utiloquity). Our sensitivity analysis demonstrates that the society should put its priority on helping eliminate poverty before investing on capacity expansions in charity organizations like food-banks. This will ensure that adding more capacity to the network will always lead to a decrease in the price of equity for the food-banks. On the other hand, we observed that encouraging the society towards charitability is always beneficial for the food-banks, albeit with diminishing returns. Finally, our experiments demonstrate that reducing poverty, as the most important element in achieving higher equity, is dependent on reducing demand variability, as opposed to the average level of demand.  
\end{abstract}


\end{frontmatter}


\pagebreak

\section{Introduction}
\label{sec:Intro}

In spite of the increase in food production by means of the modern agriculture, the undernourished population of the world has had a constant rise since 2014 reaching 10.8\% of the entire world's population, amounting to 820 million people in 2021. In the US, more than 41 million people face hunger, including nearly 13 million children. On the other hand, almost 33\% of the food produced in the world annually for human consumption is wasted, which amounts to a loss of almost 1 trillion USD. Average annual food waste per person is estimated at 95-115 kg in Europe and the US and 6-11 kg in Africa and Asia \citep{UNFAO, USDA}.

In some countries, non-profit organizations fight hunger by collecting and distributing the otherwise wasted food to the impoverished population. For instance, Feeding America is the United States' largest non-profit hunger-relief organization, operating through a nationwide network of 200 food-banks to provide food assistance to the poverty population. The organization collects food donations from national food and grocery stores, retailers, farmers and governmental agencies and distributes them in an equitable manner to its member food-banks. The food-banks subsequently distribute the food they receive from Feeding America and other sources such as local grocers, government and individuals to charitable agencies in their service regions \citep{orgut2016modeling}. Food-banks globally redistribute nearly 2.68 million metric tons of edible surplus, feeding nearly 62.5 million hungry people. In addition to feeding the poor, by saving the food from ending up in landfills, food-banks have prevented the production of 10.54 kg of greenhouse gasses \citep{GFN}.

The next tier of the supply chain, the charitable agencies, distribute the food they receive from the food-banks to poor people in their local population. Agencies are typically run by fully volunteer staff and have a certain capacity that is expressed in terms of pounds of food they can process. This capacity is generally a function of several factors such as their volunteer manpower, storage, and loading/unloading capacity. The decision to establish an agency in an area depends mainly on the census data on the poverty level in that area. Therefore each agency faces a demand proportional to the poverty population in the area it serves. 

\cite{beamon2008performance} define two pivotal performance measures to gauge humanitarian relief operations which also apply to the case of food-banks: equity and efficiency. In this paper, we address the problem of food allocation from a food-bank to charitable agencies by developing policies that are (i) equitable: as far as possible, each agency should get its fair-share of the supply, and (ii) efficient: as far as possible, waste of food across the network should be minimized.

Unlike for-profit supply chains in which meeting demand is a crucial determinant of profitability, in non-profit operations supplies are often scarce and thereby meeting demand is out of question. Therefore, allocation of the limited supply in a fair way is the main task at hand. Fairness in food-bank operations is defined by equitable distribution of food among the agencies where each agency's fair-share is proportional to the demand it serves. In our setting, we address equity by minimizing the disparity of the shares of the total supply relative to the agency's demand (represented by fill-rate) among the agencies.

Under scarce supply, volunteer-based manpower, and tight budgets, efficiency is translated into minimization of waste in non-profit network operations. Within the context of food-bank operations, limitations typically pertain to the available supplies. Therefore, efficiency in the context of food-bank operations is explained by maximum distribution of food while avoiding waste \citep{ataseven2018examination}.

Waste can arise from all three areas in the network: upstream at the supply level, downstream in the demand level, and also at the middle in the organization level. Some of the causes of waste in the upstream level are lack of coordination in terms of mismatch in both type and quantity of food donated versus the needs of the food-bank, poor information-sharing regarding the availability time of the donations, and finally the constant variability of the quantity of donations which may mandate different levels of manpower and transportation capacity from the food-bank. Causes of waste at the demand level include mismatch between the demand preferences and food provided by the agency, overestimation and underestimation of demand across the network, and poor forecasting of the available manpower for handling the allocated food at the agencies. Finally, examples of causes of waste at the organization level include lack of prioritizing schemes for using the close-to-expiration supply, unbalanced allocation of the supply across the network with little attention to capacity and demand at the agencies, and lack of coordinated communication and information sharing with the agencies. In this paper, we address waste by maximizing efficiency in the network through maximum distribution of supply while avoiding food waste at middle and downstream levels in the network.

\begin{figure}[H]
	\centering
	\includegraphics[width=0.17\textwidth]{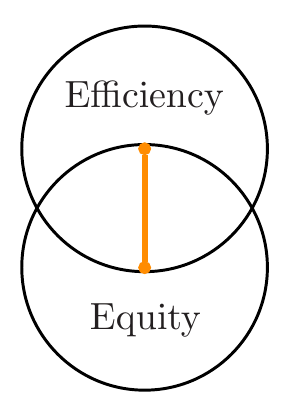}
	\caption{The dynamics of equitable and efficient policies.}
	\label{fig:Equity_Efficiency}
\end{figure}

Equity and efficiency are conflicting objectives. Figure \ref{fig:Equity_Efficiency} shows the dynamics of the two performance measures in non-profit settings. The center of each circle represents perfect accomplishment of the corresponding performance measure and the line connecting the centers represents the set of available policies that balance the non-profit organization's objectives. The choice of the policy depends on the flexibility of the organization's desired trade-off between the two performance measures, as well as practical considerations in the network. A goal of this research is to provide a flexible model to enable the food-bank decision makers to choose where in the spectrum of equity-efficiency trade-off they would like to perform.

Capacity of an agency is defined by the pounds of food it can process. Although capacity at the agencies is characterized by many factors such as volunteer staff, storage capacity, and transportation capability (among others), the number of volunteers at an agency is the key determinant of its capacity. Relying on mainly volunteer staff at the agencies, food-banks face challenges arising from human choices at the agency level. Specifically, volunteer retention and motivation are major concerns for food-banks. Although there is little control food-banks can have on volunteer behaviour, there are still strategies to ensure that food-banks can do whatever they can on their part to manage capacities in a more smart fashion. 

One of the strategies for volunteer retention that is accentuated in the literature is the avoidance of over-utilization as well as under-utilization of the capacities across the network. Over-utilization can lead to volunteer burn-out at the agencies.  Although over-utilization can increase efficiency in the short-term, it can be disadvantageous to the food-bank's capacity management in the long-run. On the other hand, under-utilization means under-engagement of the volunteers which is argued to have adverse effects on the motivation level of the volunteers, mainly arising from the effects of idleness on human psychology \citep{clary1992volunteers, agostinho2012analysis}. Therefore, another goal of this work is to address this modeling issue by creating donation allocation policies that not only address the requirements of efficiency and equity, but also implicitly consider utilization equity across the network in order to ensure volunteer retention and motivation. In the rest of the paper, we use the term ``utiloquity'' to refer to the degree by which a solution is equitable in utilizing the capacities across the agencies.

Demand variability is inherent in food-bank operations. Although a certain portion of the demand can be predicted by estimating the under-poverty population of the area served by the agencies, the walk-in behaviour of the food-banks makes them an easy food-stop for local population which in turn causes considerable variability in demand behaviour. Additionally, capacity can also be estimated at the agencies by the historical data and the fact that volunteers generally sign up before coming to an agency. Although in this paper, we assume demand and capacities at the agencies are deterministic and perfectly known when the decisions are made at the food-bank level, the third goal of this research is to design models that are robust against the variability in the data. 

Our work contributes to the literature in deterministic decision-making within the context of food-bank operations at the tactical level.
Specifically, we, (i) develop a flexible and robust model that explicitly accounts for both efficiency and equity in the objective function while implicitly optimizing utiloquity across the network, (ii) derive closed-form solutions for perfect equity and perfect efficiency, and through a numerical study show that our model dominates a benchmark from the literature in terms of efficiency, equity, and utiloquity across the full spectrum of model parameters, (iii) give the tightest spectrum of policies for food-bank decision makers to operate in, and (iv) derive managerial insights with regards to society's charitability, wealth disparity, and food-bank volunteer levels which enables managers at the food-banks to make informed decisions on their operations. 

The remainder of this paper is organized as follows: Section \ref{sec:LitReview} reviews the research most related to our paper. Sections \ref{sec:Model} and \ref{sec:ModAnalysis} introduce and analyze our model for the problem, respectively. Section \ref{sec:NumStudy} demonstrates the performance of our model compared to a benchmark from the literature and highlights our insights on the decision-making of the food-bank. Finally, Section \ref{sec:Conclusions} concludes the paper with critical insights from the study and future extensions of our work.
\section{Relevant Literature}
\label{sec:LitReview}

In contrast to for-profit organizations, the performance measures that non-profit firms seek are not solely efficiency-based. Particularly, the need for understanding the trade-off between efficiency and equity has been emphasized in the literature \citep{savas1978equity}. A stream of research to address this need in the literature has been addressed in the context of humanitarian and disaster relief operations \citep{orgut2016achieving, balcik2008last, huang2012models, taskin2010inventory, park2020supply}. A majority of research published in the context of humanitarian operations addresses the objective of efficiently distributing a resource in an equitable manner; however, most of this work only considers a disaster-related decision and not day-to-day operations as we described in our paper.

Food-banks attract a major portion of research in non-profit literature. Some of the works in food-bank research focus on operational level decisions \citep{mohan2013improving,biswal2018warehouse}. Specifically, \cite{mohan2013improving} develop a warehouse simulation model in partnership with a local food-bank and through proper demand planning, supply coordination, and logistics integration eliminate the need for extra warehouse space for handling extra volumes of supply. \cite{biswal2018warehouse} examine the benefits of RFID in terms of the effects of available rate of ordering and shrinkage recovery rate on overall costs. Our paper differs from this stream in that our work focuses on tactical level decisions of the food-bank. 

On the tactical level, two particular sub-problems have received major attention in food-bank research, i.e., (i) routing and (ii) allocation. The stream of research that focuses on routing includes scheduling of either pick-ups (from the donors) or drop-offs (at the agencies) or both, while optimizing criteria like travel distance and freshness. Some of the works in this stream are \cite{davis2014scheduling, solak2014stop, balcik2014multi, eisenhandler2019humanitarian, nair2017fair}. \cite{davis2014scheduling} model the collection and distribution of the donated supply as a set-covering problem. They assign the agencies in the network to a collection of food delivery points and show that their model reduces the food access inequity while meeting constraints such as food safety, collection frequency, and fleet capacity. Similar to \cite{davis2014scheduling}, \cite{solak2014stop} solve a food donation problem in which a central warehouse delivers food to multiple delivery sites from where the partner agencies pick up their share. Their solution involves jointly selecting a set of delivery sites, assigning agencies to these sites, and scheduling routes for the delivery vehicles. \cite{balcik2014multi} optimize the allocated food to the agencies while considering equity in the allocation and minimizing food-waste in distribution. They offer a heuristic solution to solve the problem and test the performance of their model under demand variability and changes in the supply. \cite{eisenhandler2019humanitarian} study the collection and distribution of donated supply with limited truck capacity under a maximum travel time constraint. They emphasize the importance of explicitly accounting for efficiency (maximum distribution of food) and equity (fair-shares received by the agencies) in food-bank research. \cite{nair2017fair} is another work that incorporates equity in allocation of the collected food and solves the routing and allocation model under three objective types: utilitarian (efficiency-based), egalitarian (equity-based), and deviation-based. All of the above-mentioned works study the food-bank problem with a focus on the upstream issues related to scheduling and routing. Additionally, the nature of the objective function in our setting is a fill-rate based one which accounts for all three goals of efficiency, equity, and utiloquity while these works focus on scheduling-based objective functions. 

Equity is a major consideration for many supply chains and its position has been accentuated in the literature \citep{dos2020equitable, rea2021unequal}. Within the context of food-bank research, \cite{orgut2016modeling} scrutinize the equity-efficiency trade-off, with a thorough analysis of the importance of equity in food-bank decision-making. While our work follows the same path, they account for efficiency through minimization of leftover food and for equity through maximum deviation from a perfectly equitable allocation. Our results in this model take \cite{orgut2016modeling} as a benchmark and demonstrate the dominance of our solutions to that of \cite{orgut2016modeling} in terms of efficiency, equity, and utiloquity. \cite{islam2021modeling} minimize the costs of food processing and waste at the food-bank and, similar to \cite{orgut2016modeling}, account for equity in the model through a maximum deviation tolerance from perfect equity. \cite{sengul2017modeling} focus on the effects of changes in the capacities in their previous model in \cite{orgut2016modeling} and conclude that the structure of the solution to the problem has a newsvendor behaviour. Similar to \cite{orgut2016modeling} and \cite{sengul2017modeling}, \cite{orgut2018robust} consider equity and efficiency as major elements of their model, and divide their analysis into two types of perfect equity and another that accounts for a maximum allowed deviations from perfect equity. Similar to our work in terms of consideration of a fill-rate based objective function, \cite{lien2014sequential} solve a sequential fair allocation problem, but apply a maxi-min approach to raise the minimum fill-rate across the network. Their results focus on developing near-optimal solutions for the sequence of agency visits. \cite{fianu2018markov} consider fair allocation of unknown supply among the agencies of a food-bank and identify optimal allocation rules under varying supply and demand scenarios. Similar to \cite{fianu2018markov}, \cite{alkaabneh2020unified} develop a dynamic programming approach to the food allocation problem by a food-bank to the agencies in its network and explicitly account for nutritional value of the food delivered to each agency as well as the utility and equity among the served agencies. Fairness in allocation of supply to demands across the network is studied by \cite{spiliotopoulou2021fairness}, whose numerical results suggest that when a distribution network is supply-constrained using a fill-rate based measure of equity is ideal. \cite{hynninen2020operationalization} suggest that in a resource allocation problem, cost of equity is an important measure to analyze within the contexts of utilitarian (efficiency-based) and egalitarian (equity-based) objectives. 

In this research, we focus on allocation of the available supply from a single food-bank to the agencies in its network while explicitly accounting for efficiency and equity and implicitly considering utiloquity. We provide closed-form solutions to the cases of perfect equity and perfect efficiency and provide numerical results to test the performance of our model while deriving managerial insights. 

\section{Model Formulation}
\label{sec:Model}

We consider the deterministic problem faced by a food-bank that sources $n$ charitable agencies in its network. The food-bank receives an amount of donations ($S$) which needs to be allocated to its network of agencies, where agency $i$ has a capacity level ($C_i$) to handle and distribute their share of the allocated donation among their demand ($D_i$). The aim of the food-bank is to allocate the available supply among the agencies in the most equitable and efficient manner. 

\begin{figure}[ht]
	\centering
	\includegraphics[width=0.45\textwidth]{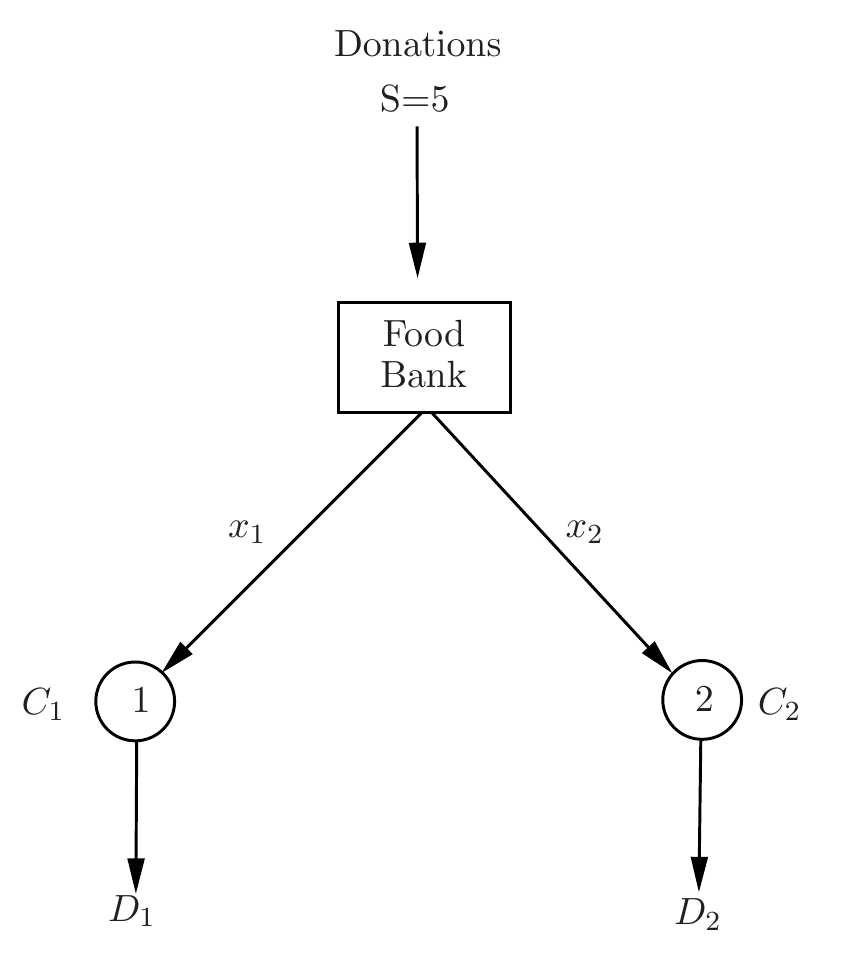}
	\caption{An example of the food-bank network with 5 units of donations and two agencies.}
	\label{fig:Setting_Example}
\end{figure}

In order to demonstrate the interaction of the two performance measures (efficiency and equity), consider a simple example presented in Figure \ref{fig:Setting_Example} in which the food-bank has received a total donation amount of 5 units and intends to allocate it to its two agencies equitably and efficiently.

Let's assume the combination of capacity pair $(C_1=1,C_2=5)$ and demand pair $(D_1=4,D_2=8)$; $x_1$ and $x_2$ are the allocated portions of the supply to the first and the second agency, respectively. Allocations to each agency cannot exceed its capacity or demand, for instance, $x_1$ cannot exceed $\min \{C_1, D_1\}=1$. We call this limit the effective demand. Let's assume the amount allocated to the first agency equals its effective demand, that is $x_1=1$. If the food-bank's policy is solely on being perfectly equitable, then one may use the demand in each area to represent its fair-share and fill the demand at exactly the same level, i.e., the fill-rates of two agencies should be the same. Therefore, $x_1/D_1=x_2/D_2=25\%$ results in $x_1=1$ and $x_2=2$. This result is perfectly equitable but not efficient. Across the network, only $(x_1+x_2)/S=60\%$ of the donations are consumed and $40\%$ are wasted, while simultaneously $75\%$ of the demand remains unsatisfied.

Now assume the combination of capacity pair $(C_1=5,C_2=3)$ and demand pair $(D_1=5,D_2=10)$. If the food-bank's policy is solely based on efficiency, then one may choose $(x_1=5,x_2=0)$ or $(x_1=2,x_2=3)$. However, neither of the two is perfectly equitable, because the first policy yields $(x_1/D_1, x_2/D_2)=(100\%, 0\%)$ and the second yields $(x_1/D_1, x_2/D_2)=(40\%, 30\%)$.

Before proceeding to the model, we present the formal structure of the supply chain under consideration in Figure \ref{fig:Setting}:

\begin{figure}[H]
	\centering
	\includegraphics[width=\textwidth]{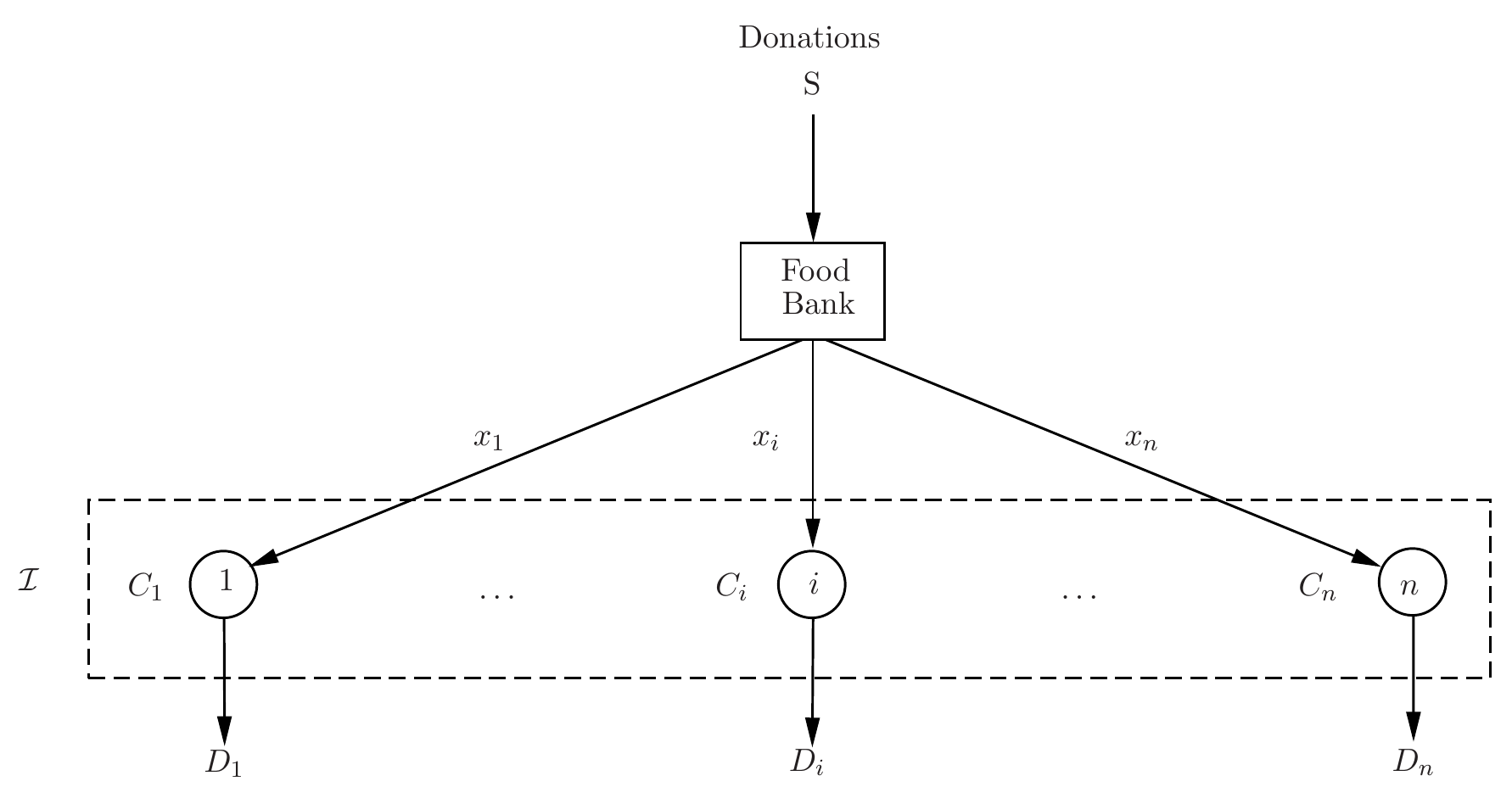}
	\caption{The overall problem setting.}
	\label{fig:Setting}
\end{figure}

In Figure \ref{fig:Setting},  $S$, $C_i$, and $D_i$ are supply at the food-bank, capacity at agency $i$, and demand at agency $i$, respectively. The decision variable $x_i$ represents the amount of donations allocated to agency $i$. These three parameters along with the decision variable $x_i$ are calculated in pounds of food. Further mathematical notations are defined in Table \ref{tab:param}.

\renewcommand{\arraystretch}{1}
\begin{table}[H]
	\centering
	\caption{Summary of the mathematical notation.}
	\begin{tabular}{llll}
		\hline
		\textbf{Notation} &  \textbf{Definition} & \textbf{Type} \\
		\hline
		$n$ & number of agencies. &  Parameter\\
		$\mathcal{I}$ & set of agencies, i.e, $\mathcal{I}=\{1, \cdots, n\}$.	&  Set\\		
		$S$ & supply of donations arrived at the food-bank. & Parameter \\
		$C_i$ & capacity at agency $i$.	& 	Parameter\\	
		$D_i$ & demand at agency $i$.& 	Parameter \\
		$\tilde{D}_i$ & effective demand at agency $i$, $\tilde{D}_i = \min \{C_i, \; D_i\}, \forall i\in \mathcal{I}$.& 	Expression \\
		$D$ & the total demand across all agencies, $D = \sum_{i \in \mathcal{I}} D_i$.& 	Expression \\	
		$x_i$ & amount of donations allocated to agency $i$.	&	Decision Variable\\
		$u_i$ & utilization rate at agency $i$, $u_i = \frac{x_i}{C_i}, \forall i\in \mathcal{I}$. & 	Expression\\
		$\beta_i$ & fill-rate at agency $i$, $\beta_i = \frac{x_i}{D_i}, \forall i\in \mathcal{I}$.& Expression\\	
		$\beta$ & maximum fill-rate among all agencies, $\beta=\max_{i\in \mathcal{I}} \beta_i$. & Decision Variable \\	
		$z_i$ & deviation from the maximum fill-rate at agency $i$, $z_i=\beta-\beta_i$. & 	Decision Variable\\
		$\theta$ & penalty for deviations from the maximum fill-rate. & 	Parameter\\
		$\theta^U$ & the tightest upper-bound of $\theta$. & 	Parameter\\	
		$\theta^L$ & the tightest lower-bound of $\theta$. & 	Parameter\\	
		\hline
	\end{tabular}
	\label{tab:param}
\end{table}

We present the mathematical programming model to the problem (Model \ref{eq:ObjFnc}) as follows

\begin{align}
Max \hspace{3pt}\sum_{i\in\mathcal{I}}(\beta_i-\theta z_i)  \label{eq:ObjFnc}
\end{align}\setcounter{equation}{0}
\begin{subnumcases}{s.t.}
   \sum_{i\in \mathcal{I}} x_i \leq S & \label{eq:Const1}\\
   x_i\leq C_i & $\forall i\in \mathcal{I}$ \label{eq:Const2}\\
   \beta-\beta_i=z_i  & $\forall i\in \mathcal{I}$ \label{eq:Const3} \\
   \beta_i = \frac{x_i}{D_i}  & $\forall i\in \mathcal{I}$ \label{eq:Const4} \\
   0\leq\beta\leq 1 & \label{eq:Const5}\\
   x_i,z_i\geq 0 & $\forall i\in \mathcal{I}$, \label{eq:Const6}
\end{subnumcases}
where $\beta_i$ is the fill-rate at agency $i$ and $\beta$ represents the maximum fill-rate among all agencies. We define fill-rate as the ratio of the allocated supply to the demand at each agency \citep{zipkin2000foundations}.

In Model \ref{eq:ObjFnc}, the objective function in Equation \eqref{eq:ObjFnc} maximizes fill-rates across the network (efficiency) while penalizing the deviations from the maximum fill-rate by parameter $\theta$ (equity). Constraint \eqref{eq:Const1} states that the total allocation across the agencies cannot exceed the total supply. Constraint \eqref{eq:Const2} sets the upper-bound on the allocation to each agency as the capacity of that agency to process food. The deviations from the maximum fill-rate for each agency is set by constraint \eqref{eq:Const3}. Constraint \eqref{eq:Const4} defines the fill-rate for each agency. Finally, constraints \eqref{eq:Const5} and \eqref{eq:Const6} set the upper-bound and lower-bound for each decision variable in Model \ref{eq:ObjFnc}. 

Using $\beta_i$ facilitates accounting for both efficiency and equity in our model. Specifically, since $\beta_i = \frac{x_i}{D_i}$, maximizing $\beta_i$ leads to maximizing $x_i$ and thereby minimizing the waste across the network in terms of food waste both at the food-bank and at the agency levels. Additionally, when all the agencies have exactly the same fill-rate such that $\beta_1=\beta_2=\cdots=\beta_n$, the allocation is perfectly equitable. As the disparity among $\beta_i$'s increases, the system becomes more and more inequitable. Therefore, $\theta$ is introduced in Model \ref{eq:ObjFnc}, indicating the degree of penalty for the deviation $\beta-\beta_i$. 

Our paper focuses on quantifying the equity-efficiency trade-off within the context of food-bank operations \citep{bertsimas2011price, bertsimas2012efficiency}. In our setting, $\theta=0$ indicates that the sole priority of the firm is efficiency, i.e., the center of the efficiency circle in Figure \ref{fig:Equity_Efficiency}. On the other hand, when $\theta$ is the arbitrarily large value $M$, it means that the sole priority of the firm is equity, i.e., the center of equity circle in Figure \ref{fig:Equity_Efficiency}. Finally, when $0<\theta<M$, it means that the firm chooses to operate with a combined approach towards efficiency and equity, i.e., the line connecting the centers of the two circles in Figure \ref{fig:Equity_Efficiency}. The values $0$ and $M$ are chosen hypothetically at this point. The tightest values for the range of $\theta$ as $\theta^L$ and $\theta^U$, which will replace $0$ and $M$, respectively, will be discussed in Section \ref{sec:ModAnalysis}. 
\section{Model Analysis}
\label{sec:ModAnalysis}
In this section, we discuss some properties of our model. The goal of our analysis is to develop the tightest upper and lower-bounds for $\theta$ in Model \ref{eq:ObjFnc} which are denoted by $\theta^U$ and $\theta^L$ respectively. We use bold face letters to indicate vectors. Additionally, all missing proofs for our theorems can be found in the appendices. 

Parameter $\theta$ in Model \ref{eq:ObjFnc} gives the degree of penalty applied to the deviation $\beta-\beta_i$, where $\beta_i$ is the fill-rate of agency $i$ and $\beta$ represents the maximum fill-rate across all the agencies. When $\theta=0$, the objective function of Model \ref{eq:ObjFnc} degenerates to $\sum_{i\in\mathcal{I}} \beta_i$. This means that the first priority of the firm is efficiency. Under such a condition the model will try to maximize the total food allocation, because the more food distributed to the agencies results in higher value for $\beta_i$'s. On the other hand, when $\theta$ is the arbitrarily large value $M$, the objective function of Model \ref{eq:ObjFnc} changes to $\sum_{i\in\mathcal{I}} \beta_i-M\sum_{i\in\mathcal{I}} (\beta-\beta_i)$. Since $M$ is arbitrarily large, obviously the second term $M\sum_{i\in\mathcal{I}} (\beta-\beta_i)$ is greater than the first term $\sum_{i\in\mathcal{I}} \beta_i \leq n$, as long as there exists a $\beta_i$ strictly less than $\beta$. This implies that the objective function value is negative if for any $i$ the inequality $\beta_i<\beta$ holds. Therefore the arbitrarily large value $M$ will force the optimization problem to find an optimal solution such that $\beta_1=\beta_2=\cdots=\beta_n=\beta$, which eliminates the effect of $M$.

To summarize, we have $\theta \in [0, M]$, where $0$ and $M$ represent perfectly efficient and perfectly equitable policies, respectively. However, in real application, the food-bank administrators may have concerns on using $M$. Particularly, the concern arises with respect to the arbitrary nature of $M$. In what follows, we give a closed-form solution for finding $\theta^U$ which conveniently enables the food-bank administrators to set the upper-bound to this value instead of $M$. In fact, we rigorously show that the value of $\theta^U$ is at most $n-1$, that is one less than the number of agencies. Next, we discuss the related definitions which culminate in the closed-form solution for $\theta^U$ in Theorem \ref{tightest bound theorem}.

\begin{defn}
\label{dfn: effective demand}
Let $\tilde{D}_i = \min \{C_i, \; D_i\}$ and $\tilde{\beta}_i = \tilde{D}_i/D_i$. $\tilde{D}_i$ and $\tilde{\beta}_i$ are called the effective demand and the maximum possible fill-rate of agency $i$ respectively.
\end{defn}

By Definition \ref{dfn: effective demand}, it is convenient to combine constraints (\ref{eq:Const2}) and (\ref{eq:Const4}) and exclude $x_i$ from our model. Model \ref{eq:ObjFnc} can be rewritten as the following concise form in which $\beta_i$ and $\beta$ are the only decision variables. Our analysis will use this concise form hereafter. What's more, when there is no confusion, $\bb=(\beta_1, \cdots, \beta_n)$ and $(\bb, \beta)$ are interchangeably called the solution to Model \ref{eq:ObjFnc-efct} because $\beta=\max_{i\in \mI} \beta_i$ which can be calculated using $\bb$.

\begin{align}
Max \hspace{3pt}\sum_{i\in\mathcal{I}}[\beta_i - \theta(\beta - \beta_i)]  \label{eq:ObjFnc-efct}
\end{align}\setcounter{equation}{1}
\begin{subnumcases}{s.t.}
   \sum_{i\in \mathcal{I}} \beta_i D_i \leq S & \label{eq:Const1-efct}\\
   \beta_i\leq \tilde{\beta}_i & $\forall i\in \mathcal{I}$ \label{eq:Const2-efct}\\
   \beta_i \leq \beta  & $\forall i\in \mathcal{I}$ \label{eq:Const3-efct} \\
   \beta\leq 1 & \label{eq:Const5-efct}
\end{subnumcases}

\begin{defn}
Given a feasible solution $\boldsymbol{\beta}=(\beta_1, \cdots, \beta_n)$, $\beta_i$ is called a binding variable if $\beta_i = \tilde{\beta}_i$.
\end{defn}

\begin{defn}
Let $\boldsymbol{\beta}^{EQ}=(\beta_1^{EQ}, \cdots, \beta_n^{EQ})$ be the optimal solution of Model \ref{eq:ObjFnc-efct} under the criterion of perfect equity, such that $\beta^{EQ}_i=\beta^{EQ}$ $(\forall i \in \mathcal{I})$ and $\sum_{i \in \mathcal{I}} \beta^{EQ}_i = n\beta^{EQ} \geq \sum_{i \in \mathcal{I}} \beta_i - M\sum_{i \in \mathcal{I}} (\beta-\beta_i)$ for any feasible solution $(\beta_1, \cdots, \beta_n)$ and an arbitrarily large value $M$. $\beta^{EQ}$ is called the perfectly equitable fill-rate.
\end{defn}

Based on the definition of $\boldsymbol{\beta}^{EQ}$, now we can rigorously define the upper-bound of $\theta$ as follows.

\begin{defn}
\label{defn: tightestbd}
$\theta^U$ is called the tightest upper-bound of $\theta$ if both of the following conditions hold:

\begin{enumerate}[(a)]
\item For any $\theta \geq \theta^U$, $\boldsymbol{\beta}^{EQ}$ is the optimal solution.
\item For any $\theta < \theta^U$, there exists a feasible solution $\boldsymbol{\beta}$ dominating $\boldsymbol{\beta}^{EQ}$ such that\\
$\sum_{i \in \mathcal{I}} \beta^{EQ}_i = n\beta^{EQ} < \sum_{i \in \mathcal{I}} \beta_i - \theta \sum_{i \in \mathcal{I}} (\beta-\beta_i)$
\end{enumerate}
\end{defn}

Next, we discuss a special case where the solution $\boldsymbol{\beta}^{EQ}$ can be found superior in both equity and efficiency. We call the problem setting under such case as utopian, rigorously defined in Definition \ref{def: utopian}.

\begin{defn}
\label{def: utopian}
Given the optimal solution $\boldsymbol{\beta}^{EQ}$ under the criterion of perfect equity, $\boldsymbol{\beta}^{EQ}$ is called utopian if the total supply is depleted, i.e.,  $\sum_{i \in \mathcal{I}} x_i = \beta^{EQ} D = S$.
\end{defn}

When supply is completely exhausted, food waste is zero. This implies that efficiency is at its maximum ($100\%$). When simultaneously all the agencies have the same fill-rate, perfect equity is also achieved. Therefore, we call this kind of optimal solution utopian, indicating perfection of both performance measures. This case is not interesting, because: (i) it is hardly true in real life.  Although our model assumes that information on capacity and demand at the agencies is perfectly known by the decision-maker prior to allocation of the supplies, capacity is variable in day-to-day operation of the food-banks, depending on the performance of the food-bank in recruitment, placement, and retention of volunteers (\citenl{clary1992volunteers}; \citenl{agostinho2012analysis}). Additionally, in most situations, the capacity in an agency is much less than the demand, resulting in a very small fill-rate. (ii) given a problem, it is straightforward to check whether the problem is utopian. If a problem is utopian, there is no need to solve the linear programming model since we give the closed-form solution of $\boldsymbol{\beta}^{EQ}$ and the condition under which it is utopian in Theorem \ref{closed form of betaEQ}.

Theorem \ref{tightest bound theorem} gives the closed-form of $\theta^U$ when the problem is not utopian:

\begin{theorem}
\label{tightest bound theorem}
Assume $\boldsymbol{\beta}^{EQ}$ is the optimal non-utopian solution of Model \ref{eq:ObjFnc-efct} under the criterion of perfect equity. Letting $m^U$ denote the number of unbinding variables in $\{\beta_1^{EQ}, \cdots, \beta_n^{EQ}\}$, $\theta^U$ is calculated as follows.
\begin{align}
\theta^U = \frac{m^U}{n-m^U}
\end{align}
\end{theorem}

\noindent \underline{\emph{Proof}}

We prove Theorem \ref{tightest bound theorem} by parts. If $m^U=0$, it means that all the decision variables are binding such that $\tilde{\beta}_i = \beta^{EQ}_i = \beta^{EQ}$. On the other hand, by the constraints (\ref{eq:Const2-efct}) and (\ref{eq:Const3-efct}), the objective function value of Model \ref{eq:ObjFnc-efct} is $\sum_{i\in\mathcal{I}}[\beta_i - \theta(\beta - \beta_i)] \leq \sum_{i\in\mathcal{I}} \beta_i \leq \sum_{i\in\mathcal{I}} \tilde{\beta}_i = \sum_{i\in\mathcal{I}} \beta_i^{EQ}$. Hence, $\boldsymbol{\beta}^{EQ}$ is the optimal solution for any $\theta \geq \theta^U = \frac{m^U}{n-m^U} = 0$.

Now, let's consider $m^U>0$. We prove this by the definition of $\theta^U$ (Definition \ref{defn: tightestbd}). Notice that $\boldsymbol{\beta}^{EQ}$ is not utopian, so $\sum_{i\in \mathcal{I}} \beta_i^{EQ} D_i < S$.
\begin{enumerate}[i.]
\item Given a $\theta<\frac{m^U}{n-m^U}$, there exists a $\boldsymbol{\beta} \neq \boldsymbol{\beta}^{EQ}$ that dominates $\boldsymbol{\beta}^{EQ}$.\\
Denote by $\mathcal{A} = \{i \mid \tilde{\beta}_i = \beta^{EQ} \}$, the set of binding indices. We can construct a solution $\boldsymbol{\beta}$ by a sufficiently small variable $\alpha>0$ such that $\beta_i = \beta^{EQ}$ for any $i \in \mathcal{A}$, $\beta_i = \beta^{EQ} + \alpha$ for any $i \in \mathcal{I}/\mathcal{A}$, and $\beta=\beta^{EQ}+\alpha$. Since $\alpha$ is sufficiently small, it is trivial that ($\boldsymbol{\beta}$, $\beta$) is a feasible solution of Model \ref{eq:ObjFnc-efct}. On the other hand, letting $|\cdot|$ be the cardinality of a set, we have
\begin{align*}
 &\sum_{i\in\mathcal{I}}[\beta_i - \theta(\beta - \beta_i)] -  \sum_{i\in\mathcal{I}} \beta^{EQ}_i\\
=&\sum_{i\in\mathcal{I}}[\beta_i - \theta(\beta - \beta_i) -  \beta^{EQ}] \\
=&\sum_{i\in\mathcal{I}}[\beta_i - \beta^{EQ}] - \sum_{i\in\mathcal{I}} \theta(\beta - \beta_i) \\
=&|\mathcal{I}/\mathcal{A}| \alpha - |\mathcal{A}| \theta \alpha  \\
=&[m^U- \theta (n-m^U)] \alpha \\
>&[m^U- (n-m^U) \frac{m^U}{n-m^U}] \alpha = 0
\end{align*}

Hence, we found a feasible solution $\boldsymbol{\beta}$ dominating $\boldsymbol{\beta}^{EQ}$, for any $\theta <\frac{m^U}{n-m^U}$.

\item Given a $\theta \geq \frac{m^U}{n-m^U}$, $\boldsymbol{\beta}^{EQ}$ is the optimal solution.\\
Let ($\boldsymbol{\beta}$, $\beta$) be any feasible solution. Because $n-m^U$ variables are binding and $\beta_i \leq \tilde{\beta}_i = \beta^{EQ}$ for the binding variables (in set $\mathcal{A})$, we have $\sum_{i\in\mathcal{I}} \theta (\beta-\beta_i) \geq \sum_{i\in\mathcal{A}} \theta (\beta-\beta_i) \geq (n-m^U)\theta (\beta-\beta^{EQ})$ and $\sum_{i\in\mathcal{I}} (\beta_i-\beta^{EQ}) \leq \sum_{i\in\mathcal{I}/\mathcal{A}} (\beta_i-\beta^{EQ}) \leq m^U (\beta-\beta^{EQ})$. Consider the difference between the objective function values as follows.
\begin{align*}
 &\sum_{i\in\mathcal{I}} \beta^{EQ}_i - \sum_{i\in\mathcal{I}}[\beta_i - \theta(\beta - \beta_i)]\\
=&\sum_{i\in\mathcal{I}} \theta(\beta - \beta_i) -  \sum_{i\in\mathcal{I}} (\beta_i - \beta^{EQ}) \\
\geq & (n-m^U) \theta (\beta-\beta^{EQ}) - m^U (\beta-\beta^{EQ}) \\
\geq & (n-m^U) \frac{m^U}{n-m^U}(\beta-\beta^{EQ}) - m^U (\beta-\beta^{EQ}) = 0
\end{align*}

Hence, $\boldsymbol{\beta}^{EQ}$ dominates any feasible solution $\boldsymbol{\beta}$ and is of course the optimal.

Combining the above two parts of proof, we show that $\theta^U = \frac{m^U}{n-m^U}$ is the tightest upper-bound by the definition.
\end{enumerate}

\begin{flushright}   Q.E.D.   \end{flushright}

\begin{theorem}
\label{closed form of betaEQ}
Let $\beta^a = \frac{S}{D}$, $\beta^b = \displaystyle\min_{i \in \mathcal{I}} \tilde{\beta}_i$, and $\mathcal{J} = \{j \mid \tilde{\beta}_j=\beta^b, \; \forall  j \in \mathcal{I}\}$. The following results hold.
\begin{enumerate}[i.]
\item If $\beta^a \leq \beta^b$, the problem is utopian and $\beta^{EQ}=\beta^a$.
\item If $\beta^a > \beta^b$, then $\beta^{EQ}=\beta^b$ and
\begin{align}
m^U=n-|\mathcal{J}|
\end{align}
where $|\cdot|$ represents the cardinality of a set.
\end{enumerate}
\end{theorem}

Theorem \ref{closed form of betaEQ} actually gives the closed-form of $\beta^{EQ}$.

So far, we have discussed how to find the tightest upper-bound $\theta^U$. Instead of an unknown arbitrarily large ``M'', $\theta^U$ conveniently enables the food-bank decision makers to know where exactly in the spectrum of values of $\theta$ perfect equity begins. Similar to the case of the upper-bound on $\theta$, we next discuss how the decision-maker can adopt $\theta^L$ (the lower-bound of $\theta$) to replace $0$. $\theta^L$ essentially marks where exactly in the spectrum of values of $\theta$ perfect efficiency begins. Before proceeding to the closed-form of $\theta^L$, we first characterize the nature of $\theta^L$ in Definition \ref{defn: tightestlowbd}.

So far, we have discussed how to find the tightest upper-bound $\theta^U$. Instead of an unknown arbitrarily large ``M'', $\theta^U$ conveniently enables the food-bank decision makers to adopt polices that contribute to increased equity. Similar to the case of the upper-bound on $\theta$, we next discuss how the decision-maker can  adopt $\theta^L$ (the lower-bound of $\theta$) to replace $0$. Before proceeding to the closed-form of $\theta^L$, we first characterize the nature of $\theta^L$ in Definition \ref{defn: tightestlowbd}.

\begin{defn}
\label{defn: tightestlowbd}
Assuming $\xi$ is the optimal objective function value when $\theta=0$, $\theta^L$ is called the tightest lower-bound of $\theta$ if both of the following conditions hold:
\begin{enumerate}[(a)]
\item For any $\theta \leq \theta^L$, if $\boldsymbol{\beta}=(\beta_1, \cdots, \beta_n)$ is the optimal solution, then $\sum_{i \in \mathcal{I} } \beta_i = \xi$.
\item For any $\theta > \theta^L$, if $\boldsymbol{\beta}=(\beta_1, \cdots, \beta_n)$ is the optimal solution, then $\sum_{i \in \mathcal{I} } \beta_i < \xi$.
\end{enumerate}
\end{defn}

The reason that we do not specify the optimal solution for $\theta=0$ is because there may exist multiple or infinitely many optimal solutions when efficiency is the sole criterion in the objective function. For instance, let $S=C_i = 50$ for any $i \in \mathcal{I}=\{1,2,3\}$, $D_1=20$, and $D_2=D_3=50$. When $\theta=0$ the optimal solution is $(\beta_1^*, \beta_2^*, \beta_3^*) = (100\%, \; x/50, \; 60\%-x/50)$. $x$ can be any value between 0 and 30 so there exists infinitely many optimal solutions. However, the optimal objective function value $\xi$ is unique. This is why we use $\xi$ to define the tightest lower-bound. 

The unique optimal solution $\boldsymbol{\beta}^{EQ}$ provided the clue to find the tightest upper-bound. Similarly, in what follows we show that although the problem may not have a unique optimal solution for $\theta=0$, there exists a dominant solution which facilitates our analysis. Specifically, we develop Algorithm \ref{alg:domiEfficient} to find this unique optimal dominant solution. Intuitively, the rationale of the algorithm is motivated by that of the greedy method. Specifically, when $\theta=0$, we only focus on the efficiency, that is, maximizing the fill-rates $\beta_i$'s. Considering $\beta_i = x_i/D_i$, the agency with the smallest demand $D_i$ should be filled first since the same amount of food supply generates a higher fill-rate. So we divide agency set $\mathcal{I}$ into subsets denoted by $\mathcal{I}_i$'s. In each subset $\mathcal{I}_i$, all the agencies have the same demand such that $D_{k_1}=D_{k_2}=D_i$ for any $k_1, \, k_2 \in \mathcal{I}_i$. Assuming there are $\tilde{n}$ different subsets, we have $\mathcal{I}=\bigcup_{i=1}^{\tilde{n}} \mathcal{I}_i$, which reduces $n$ agencies to $\tilde{n}$ groups. Without loss of generality, assume $D_1 \leq D_2 \leq \cdots \leq D_{\tilde{n}}$. Our algorithm attempts to fill the first group (subset) of agencies, then the second, ..., until there is no food supply is left on-hand. Next task is to characterize the order by which the agencies in the same group $\mathcal{I}_i$ are filled. Although the agencies in the same group have identical demands, they differ in effective demand which is the minimum value of demand and capacity for each agency. Therefore, we can further divide $\mathcal{I}_i$ into subsets $\mathcal{I}_{ij}$. In each subset $\mathcal{I}_{ij}$, all the agencies have not only the same demand but also the same effective demand such that $\tilde{D}_{k_1}=\tilde{D}_{k_2}=\tilde{D}_{ij}$ for any $k_1, \, k_2 \in \mathcal{I}_{ij}$. Without loss of generality, assume there are $\tilde{n}_i$ different subsets in group $\mathcal{I}_i$ and $\tilde{D}_{i1} \leq \tilde{D}_{i2} \leq \cdots \leq \tilde{D}_{i\tilde{n}_i}$. When filling the agencies in the group $\mathcal{I}_i$, however, we do not fill $\mathcal{I}_{ij}$ one by one in ascending order. This is because even if two agencies belong to different subsets (i.e., different effective demands) in the group $\mathcal{I}_i$ they still have the same demand. Hence, allocating all the supply exclusively to one subset $\mathcal{I}_{ij}$ is identical to spreading the supply over different $\mathcal{I}_{ij}$'s simultaneously. However, the second approach is preferred as it favors equity in addition to efficiency. This will benefit the solution when the penalty factor $\theta$ becomes positive to penalize the inequity. This procedure is formally delineated in Algorithm \ref{alg:domiEfficient}. However, the summary of how Algorithm \ref{alg:domiEfficient} fills the subset $\mathcal{I}_i$ is given in the following three steps:
\begin{enumerate}[(i)]
\item Initialize $j=1$
\item Denote by $\mathcal{I}'_j = \bigcup_{k=j}^{\tilde{n}_i} \mathcal{I}_{ij}$ the subsets not yet fully filled. Check whether the remaining supply can increase the fill-rate of all the agencies in $\mathcal{I}'_j$ with the same amount $\tilde{\beta}_{ij}$. If yes, allocate food to all the agencies in $\mathcal{I}'_j$ to raise the fill-rates by $\tilde{\beta}_{ij}$. If no, then distribute the remaining supply and equally raise the fill-rates for all the agencies in $\mathcal{I}'_j$.
\item Adjust $S$ and $\tilde{\beta}_{ij}$. If the remaining supply $S$ reaches $0$, stop; otherwise, $j=j+1$ and go to step (ii).
\end{enumerate}

\begin{algorithm}[ht]
\caption{The Dominant Efficient Solution}
\label{alg:domiEfficient}
{\fontsize{10}{16}\selectfont
\begin{flushleft}
\textbf{Input: } the subsets $\mathcal{I}=\bigcup_{i=1}^{\tilde{n}} \mathcal{I}_i$ and $\mathcal{I}_i=\bigcup_{j=1}^{\tilde{n}_i} \mathcal{I}_{ij}$; the parameters $S$, $D_i$, $\tilde{D}_{ij}$, and $\tilde{\beta}_{ij}$ ($i=1, \cdots, \tilde{n}$; $j=1, \cdots, \tilde{n}_i$).\\
\textbf{Output: } $\boldsymbol{\beta}^{EF}=(\beta_1^{EF}, \cdots, \beta_n^{EF})$ \\
\textbf{Initialization: } let $\beta_i^{EF}=0$ for any $i \in \mathcal{I}$.\\
\end{flushleft}
	\begin{algorithmic}[1]
	\label{Alg:DominantSol}
		\Procedure {Update}{$S$, $\mathcal{I}$, $D_i$, $\tilde{D}_{ij}$, $\tilde{\beta}_{ij}$}
			\For{$i = 1$ \textbf{to} $\tilde{n}$}
				\If {S=0}
					\State \textbf{break}
				\Else
				    \For{$j = 1$ \textbf{to} $\tilde{n}_i$}
    				    \State $\mathcal{I}'_j = \bigcup_{k=j}^{\tilde{n}_i} \mathcal{I}_{ik}$
    				    \If {$|\mathcal{I}'_j|\tilde{D}_{ij} \leq S$}
        				    \For{$\forall k \in \mathcal{I}'_j$} {$\beta_k^{EF} = \beta_k^{EF} + \tilde{\beta}_{ij}$} \EndFor
        				    \For{$k=j$ \textbf{to} $\tilde{n}_i$} {$\tilde{D}_{ik}=\tilde{D}_{ik}-\tilde{D}_{ij}$} \EndFor
        				    \State $S=S-|\mathcal{I}'_j|\tilde{D}_{ij}$
    				    \Else
        				    \For{$\forall k \in \mathcal{I}'_j$} {$\beta_k^{EF} = \beta_k^{EF} + \frac{S/|\mathcal{I}'_j|}{D_i}$} \EndFor
        				    \For{$k=j$ \textbf{to} $\tilde{n}_i$} {$\tilde{D}_{ik}=\tilde{D}_{ik}-S/|\mathcal{I}'_j|$} \EndFor
        				    \State $S=0$    				        
    				    \EndIf
    				    \If {$S=0$}
    					    \State \textbf{break}
    				    \EndIf
				    \EndFor
				\EndIf
			\EndFor
		\EndProcedure
		\Return $\boldsymbol{\beta}^{EF}$
	\end{algorithmic}
	}
\end{algorithm}

The optimal solution constructed by Algorithm \ref{alg:domiEfficient} dominates all the optimal solutions solely focusing on efficiency. In other words, if we collect all the optimal solutions for $\theta=0$ in a set $\mathcal{B}$, then this unique optimal solution dominates the rest of solutions in $\mathcal{B}$ for a positive $\theta$. This is formally stated in Theorem \ref{dominant efficient optimal}.

\begin{theorem}
\label{dominant efficient optimal}
Let $\mathcal{B}$ be the set of all optimal solutions when $\theta=0$ and $\boldsymbol{\beta}^{EF}=(\beta^{EF}_1, \cdots, \beta^{EF}_n)$ the dominant solution constructed by Algorithm \ref{alg:domiEfficient}. For any $\boldsymbol{\beta} \in \mathcal{B}$ and $\theta\geq 0$, $\boldsymbol{\beta}^{EF}$ dominates $\boldsymbol{\beta}$ such that $\sum_{i \in \mathcal{I}} [\beta^{EF}_i-\theta(\beta^{EF}-\beta^{EF}_i)] \geq \sum_{i \in \mathcal{I}} [\beta_i-\theta(\beta-\beta_i)]$.
\end{theorem}

Since the optimal solution $\bb^{EF}$  constructed by Algorithm \ref{alg:domiEfficient} dominates all the other optima for $\theta=0$, to find the tightest lower-bound, we only have to show that while $\theta$ becomes positive whether $\bb^{EF}$ is still optimal. Based on this point and Definition \ref{defn: tightestlowbd}, we present Theorem \ref{tightest lower-bound}, presenting the tightest lower-bound, the proof of which is similar to that of the tightest upper-bound.

\begin{theorem}
\label{tightest lower-bound}
Let $\bb^{EF}$ be the dominant optimal solution constructed by Algorithm \ref{alg:domiEfficient}, $\beta^{EF}=\displaystyle\max_{i \in \mI} \beta^{EF}_i$ and $\mathcal{K}=\{i \mid \beta^{EF}_i=\beta^{EF}, \; i \in \mI\}$. Let $D_m = \displaystyle\min_{i \in \mathcal{L}} D_i$, where $\mathcal{L}=\{i \mid \beta_i^{EF}<\tilde{\beta}_i, \; i \in \mI/\mathcal{K}\}$ -- the set of agencies not fully filled except those in $\mathcal{K}$. If the problem is non-utopian and supply-constrained such that $S < \displaystyle\sum_{i \in \mI} \tilde{D}_i$, the tightest lower-bound $\theta^L$ is calculated as follows.
\begin{align*}
    \theta^L = \frac{m^L}{n-m^L},
\end{align*}
where,
\begin{align*}
m^L=
    \begin{cases}
    |\mathcal{K}|-\displaystyle\sum_{i\in \mathcal{K}} D_i/D_m & \text{if $\mathcal{L}$ is not empty}\\
    |\mathcal{K}| & \text{if $\mathcal{L}$ is empty}
    \end{cases}
\end{align*}
\end{theorem}

\noindent \underline{\emph{Proof}}

Because the problem is not utopian, it is easy to see $n\neq m^L$. Hence, $\theta^L$ is well defined. We prove this theorem by parts based on the Definition \ref{defn: tightestlowbd}. 
\begin{enumerate}[i.]
\item Given a $\theta>\theta^L=\frac{m^L}{n-m^L}$, there exists a feasible solution $(\bb, \beta)$ such that $\sum_{i \in \mI} \beta_i<\beta^0$ and it dominates $(\bb^{EF}, \beta^{EF})$ created by Algorithm \ref{alg:domiEfficient}.\\
Denote $\mathcal{K}=\{i \mid \beta^{EF}_i=\beta^{EF}, \; i \in \mI\}$, the set of agencies with the highest level of fill-rate; and $\mathcal{A}=\{i \mid D_i=D_m, \; i\in \mathcal{L}\}$, the set of agencies with minimum demand among those can be further filled. Assume $\mathcal{L}$ is not empty. Denoting $\rho=\frac{\sum_{i\in \mathcal{K}}D_i}{D_m}$, we construct a solution $\boldsymbol{\beta}$ by a sufficiently small $\alpha>0$ such that $\beta_i = \beta^{EF}_i-\alpha$ for any $i \in \mathcal{K}$, $\beta_i=\beta^{EF}_i+\alpha\rho/|\mathcal{A}|$ for any $i \in \mathcal{A}$ and $\beta_i=\beta_i^{EF}$ for any $i \in \mI/(\mathcal{K}\cup \mathcal{A})$. Since $\alpha$ is sufficiently small and $S<\sum_{i\in \mathcal{I}} \tilde{D}_i$, it is trivial that ($\boldsymbol{\beta}$, $\beta$) is a feasible solution of Model \ref{eq:ObjFnc-efct}. Therefore $\beta=\beta^{EF}-\alpha$, we have
\begin{align*}
 &\sum_{i\in\mathcal{I}}[\beta_i - \theta(\beta - \beta_i)] -  \sum_{i\in\mathcal{I}}[\beta_i^{EF} - \theta(\beta^{EF} - \beta_i^{EF})]\\
=&\sum_{i\in\mathcal{I}}(\beta_i - \beta^{EF}_i) + \theta \sum_{i\in\mathcal{I}} [(\beta^{EF}-\beta) + (\beta_i - \beta_i^{EF})] \\
=&(\alpha\rho-\alpha|\mathcal{K}|)+\theta (\alpha n+\alpha\rho-\alpha|\mathcal{K}|)\\
=&-\alpha m^L+\alpha \theta (n-m^L)\\
>&-\alpha m^L+\alpha (n-m^L) \frac{m^L}{n-m^L}  = 0
\end{align*}
When $\mathcal{L}$ is empty, we can construct a solution $\boldsymbol{\beta}$ by a sufficiently small $\alpha>0$ such that $\beta_i = \beta^{EF}_i-\alpha$ for any $i \in \mathcal{K}$ and $\beta_i=\beta_i^{EF}$ for any $i \in \mI/\mathcal{K}$. The proof of this case is similar to the one when $\mathcal{L}$ is not empty. Hence, we found a feasible solution $\boldsymbol{\beta}$ dominating $\boldsymbol{\beta}^{EF}$, for any $\theta > \theta^L$. Moreover, considering the procedure how we construct $\bb^{EF}$, $\alpha|\mathcal{K}|>\alpha\rho$ because $D_i<D_m$ for any $i\in \mathcal{K}$. Therefore, $\sum_{i \in \mI} \beta_i<\beta^0=\sum_{i \in \mI} \beta_i^{EF}$. 

\item Given a $\theta\leq \theta^L=\frac{m^L}{n-m^L}$, $\bb^{EF}$ is the optimal solution.\\ 
When $\theta=0$, of course $\bb^{EF}$ is the optimal. If $0<\theta\leq \theta^L=\frac{m^L}{n-m^L}$, assume $\mathcal{L}$ is not empty and consider any feasible solution $\bb$ other than $\bb^{EF}$. Because $\sum_{i \in \mathcal{I}} \beta^{EF} \geq \sum_{i \in \mathcal{I}} \beta$ and the difference between objective function values of $\bb^{EF}$ and $\bb$ is 
$$\Delta(\bb)=(1+\theta)\sum_{i\in\mathcal{I}}(\beta_i^{EF} - \beta_i) + \theta \sum_{i\in\mathcal{I}} (\beta-\beta^{EF}),$$
if $\beta \geq \beta^{EF}$ then $\bb^{EF}$ dominates $\bb$ since $\Delta \geq 0$. Now, consider the case $\beta < \beta^{EF}$. \\ 
For a $\bb$, if $\exists i\in \mathcal{K}$ such that $\beta_i < \beta$, we can always find a solution $\bb'$ with $\beta'_i=\beta$ but dominating $\bb$, i.e., $\Delta(\bb') \leq \Delta(\bb)$. To obtain $\beta_i$, we have to first reduce $\beta_i^{EF}$ (which is equal to $\beta^{EF}$) to $\beta$ and then further reduce it to $\beta_i$. However, the further reduction does not change $\beta$ but only saved the supply by $D_i(\beta-\beta_i)$. Based on the way how Algorithm \ref{alg:domiEfficient} constructs $\bb^{EF}$, the saved supply could at most increase the total fill-rates of $\bb$ by $D_i/D_m(\beta-\beta_i)$ which is less than or equal to $\beta-\beta_i$ because $D_i \leq D_k$ for any $k \in \mathcal{L}\cup\mathcal{K}$. So the further reduction cannot reduce the value of first term of $\Delta$ and does not change the value of the second term. This implies that a $\bb'$ without further reducing the fill-rate of agency $i$ from $\beta$ to $\beta_i$ dominates $\bb$. WLOG, we assume $\beta_i=\beta, \; \forall i \in \mathcal{K}$.
\\
Let $\mathcal{K}_1=\{i \mid \beta_i<\beta^{EF}_i, \; i \in \mI\}$ and $\mathcal{K}_2=\{i \mid \beta_i>\beta^{EF}_i, \; i \in \mI\}$. Trivially, $\mathcal{K}_1 \supseteq K$ notice that $\beta<\beta^{EF}$. Furthermore, because $\sum_{i\in \mathcal{K}_1} D_i(\beta^{EF}_i-\beta_i)=\sum_{i\in \mathcal{K}_2} D_i(\beta_i-\beta^{EF}_i)$ by the assumption that the problem is supply-constrained, we observe
\begin{align*}
    \sum_{i\in \mathcal{K}_2} (\beta_i-\beta_i^{EF}) \leq \sum_{i\in \mathcal{K}_1} D_i/D_m (\beta^{EF}_i-\beta_i) =
    \rho(\beta^{EF}-\beta)+\sum_{i\in \mathcal{K}_1/\mathcal{K}} D_i/D_m (\beta^{EF}_i-\beta_i),
\end{align*}
where the second equality comes from the fact such that $\beta_i^{EF} -\beta_i = \beta^{EF}-\beta$ for any $i \in \mathcal{K}$. Hence, we have
\begin{align*}
    \Delta(\bb) &= (1+\theta)\sum_{i\in\mathcal{I}}(\beta_i^{EF} - \beta_i) + \theta \sum_{i\in\mathcal{I}} (\beta-\beta^{EF})\\
     &= (1+\theta)[|\mathcal{K}|(\beta^{EF}-\beta)+\sum_{i \in \mathcal{K}_1/\mathcal{K}} (\beta^{EF}_i-\beta_i)-\sum_{i \in \mathcal{K}_2} (\beta^{EF}_i-\beta_i)]+n\theta(\beta-\beta^{EF})\\
     &\geq (1+\theta)[(|\mathcal{K}|-\rho)(\beta^{EF}-\beta)+\sum_{i\in \mathcal{K}_1/\mathcal{K}} (1-D_i/D_m) (\beta^{EF}_i-\beta_i)]+n\theta (\beta-\beta^{EF})\\
     &\geq (1+\theta)m^L(\beta^{EF}-\beta)+n\theta (\beta-\beta^{EF})\\
     &\geq \frac{nm^L}{n-m^L}(\beta^{EF}-\beta)-\frac{nm^L}{n-m^L}(\beta^{EF}-\beta)=0
\end{align*}
Therefore, $\beta^{EF}$ is the optimal solution. For the case when $\mathcal{L}$, we have $\mathcal{K}_2=\emptyset$. The proof of this case is similar and more straightforward.
\end{enumerate}

In sum, the above two parts together complete the proof.

\begin{flushright}   Q.E.D.   \end{flushright}

Theorem \ref{tightest lower-bound} applies to the case where the total supply cannot satisfy all the agencies' effective demands. The case with sufficient supply such that $S \geq \sum_{i \in \mI} \tilde{D}_i$, is similar to the utopian problem and rarely realistic. In early 2020, the novel coronavirus (COVID-19) began to spread across the United States, and one of the consequences was an economic recession that resulted in more food insecurity – the lack of access to sufficient food because of limited financial resources. Feeding America estimates that 45 million people (1 in 7), including 15 million children (1 in 5), may have experienced food insecurity in 2020 \citel{FA}. If supply is abundant, food-bank administrators can simply satisfy all the agencies at their maximum needs.  

In summary, this section presents the tightest spectrum $[\theta^L, \theta^U]$ of policies for food-bank decision makers to operate in. We also give the closed-form solutions corresponding to $\theta^U$ and $\theta^L$. As a result, the food-bank administrators can pick a policy (defined by a $\theta \in [\theta^L, \theta^U]$) according to their organization's preference towards efficiency and the equity. Next in Section \ref{sec:NumStudy}, using real data from Feeding America, we demonstrate the performance of our model against a benchmark from the literature and provide critical insights for the food-bank decision-makers.

\section{Computational Study and Insights}
\label{sec:NumStudy}

In this section, we summarize the performance of our model compared to that of \cite{orgut2016modeling} as our benchmark from the literature. We programmed both models in Python under Google Colaboratory. We further demonstrate how changes in charitability as well as volunteer levels in the society affect the equity-efficiency trade-off for the food-banks. Specifically, we provide guidelines on the policies the society can adopt in order to facilitate the food-bank's path towards achieving higher levels of equity while sacrificing the least from their efficiency.

We use two performance measures of equity and efficiency defined as follows.

\begin{itemize}
    \item Efficiency: the ratio of the total food distributed towards satisfying the demand to the total available supply, calculated according to the following:
     \begin{equation}
       \label{eq:EFDef} \text{Efficiency}=\frac{S-(\text{Agency Waste} + \text{Food-Bank Waste})}{S}
    \end{equation}  

We define agency waste as the amount of food which has been sent to an agency, but not distributed towards satisfying the demand. We would like to point out that agency waste can happen only if the food-bank allocates food to a particular agency beyond its demand. Our model avoids agency waste via model constraints; however, it is possible to have agency waste in \cite{orgut2016modeling}. Similarly, food-bank waste is defined as the amount of food that is left at the food-bank and not allocated to the agencies.
    \item Equity: the degree to which each agency serving the impoverished population receives its fair-share of the supply $S$, calculated according to an equity measure.

In order to calculate the equity performance measure, we need to utilize a measure of equity. However, to avoid divisions by zero, in the literature, various measures of inequity have been utilized instead of equity. In this paper, we utilize some of the most widely used measures of inequity for a given vector of real values $\boldsymbol{E}=(E_1,E_2, \dots, E_n)$ having the mean of $\bar{E}$, maximum of $E_{max}$, and minimum of $E_{min}$ as follows.
\begin{align}
    &\text{Gini Coefficient} =\frac{\sum_{i\in \mathcal{I}} \sum_{j\in \mathcal{I}}\left| E_i-E_j \right|}{2n^2\bar{E}}  \label{eq:IEQ_Gini}\\
    &\text{Coefficient of Variation} =\frac{\sqrt{\frac{\sum_{i \in \mathcal{I}}(E_i - \bar{E})^2}{n}}}{\bar{E}}  \label{eq:IEQ_CV}\\  
    &\text{Variance} =\frac{\sum_{i \in \mathcal{I}}(E_i - \bar{E})^2}{n}  \label{eq:IEQ_Var}\\ 
    &\text{Mean Absolute Deviation} =\frac{\sum_{i \in \mathcal{I}}\left|E_i - \bar{E}\right|}{n}  \label{eq:IEQ_MAD}\\ 
    &\text{Range} = E_{max} - E_{min} \label{eq:IEQ_Range}
\end{align}
\end{itemize}

In this section, to calculate inequity we put our focus on Gini coefficient, the most widely used measure of inequity in social welfare \citep{marsh1994equity}. The rest of the measures used for comparing our model to that of \cite{orgut2016modeling} are reported in \ref{app: EQ_OtherMeasures}.

Gini coefficient measures inequity as the ratio of the area between the line of perfect equity and the Lorenz curve and the total area under the line of perfect equity, thereby it is a measure between zero and one \citep{gini1912variabilita}. Therefore, to calculate equity in this section, we subtract the Gini coefficient from one.

In order to compare the performance of the two models, we use the data from 2020 for a service area in North Carolina. The agencies working with Feeding America serve a single meal per day, every day throughout the month. To estimate the demand in pounds of food at each agency (represented by $D_i$ in our model), we use 1.2 lb consumption per person per serving multiplied by the number of poor persons served in the population. To calculate the capacity at each agency (represented by $C_i$ in our model), we use the $C_i/D_i$ ratios from \cite{orgut2016modeling}. The total amount of food donations to be allocated to the 34 agencies in the area under study is 2,838,584 lb, which is represented by parameter $S$ in our model. These three parameters form our three nominal parameters required to generate the test cases for our analysis.

\begin{figure}[H]
	\begin{center}
		\subfigure[\emph{low} demand variability]{\includegraphics[width=0.48\textwidth]{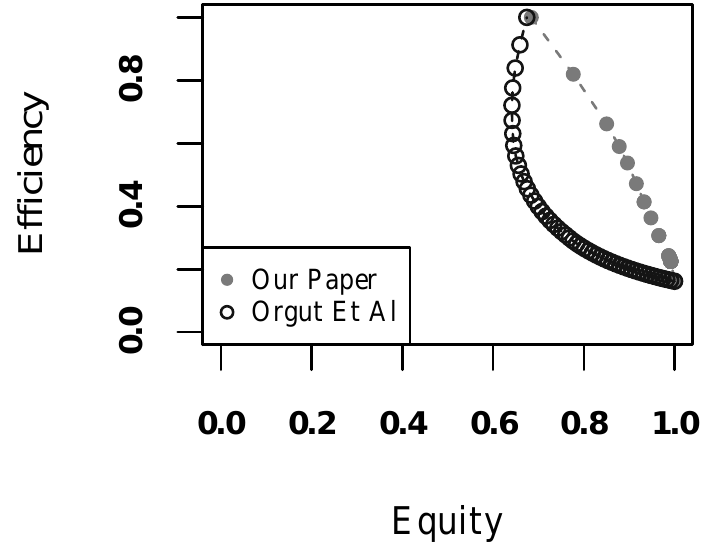} \label{fig:Bin1-EF-EQ_Gini}}
		\subfigure[\emph{high} demand variability]{\includegraphics[width=0.48\textwidth]{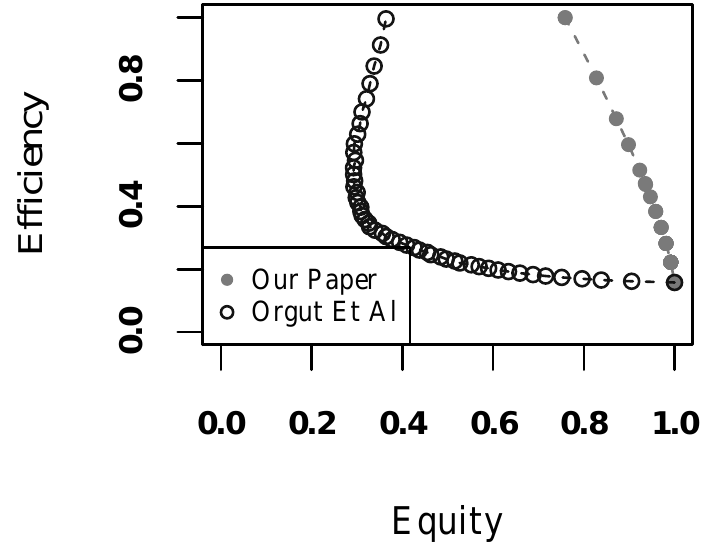} \label{fig:Bin2-EF-EQ_Gini}}
		\caption{Efficient frontier using one minus the Gini Coefficient as the equity measure.}
		\label{fig:EF-EQ_Gini}
	\end{center}
\end{figure}

To test the performance of the models against the variability in the data, using the nominal demand vector, we generate a total of 2000 new demand realizations. Specifically, we use truncated normal distribution to generate the demand vectors within two variation categories of high and low, each containing 1000 random realizations. Keeping the means the same as the nominal demand vector, we generate our low variability category to contain demand vectors with standard deviations below 20\% of that of the nominal, and our high variability category to have at least 180\% more standard deviations than that of the nominal. Additionally, we make sure that all our instances satisfy the two criteria for non-triviality discussed in Section \ref{sec:ModAnalysis}: (i) they are not utopian (from Definition \ref{def: utopian}); and (ii) they are supply-constrained (from Theorem \ref{tightest lower-bound}). Every point in Figures \ref{fig:EF-EQ_Gini}, \ref{fig:EQ-UEQ_Gini}, \ref{fig:PrEQ_S}, and \ref{fig:PrEQ_C}, along with the figures in the appendices are averaged over the 1000 demand realizations. 

Next, we determine the range of the auxiliary parameters $\theta$ (for our model) and $K$ (for \cite{orgut2016modeling} model). In detail, the lower-bound of $\theta$ in our model considers only efficiency as the goal of the objective function and it is calculated using Theorem \ref{tightest lower-bound}. On the other hand, the upper-bound of $\theta$ only takes equity into account in our objective function and is calculated based on Theorem \ref{tightest bound theorem} in our model. 

To determine the corresponding range for $K$, we note that the mathematical formulation of the problem according to \cite{orgut2016modeling} is as follows.

\begin{align}
Min \hspace{3pt}P  \label{eq:ObjFnc_Orgut}
\end{align}\setcounter{equation}{10}
\begin{subnumcases}{s.t.}
   \left| \frac{x_i}{\sum_{i\in \mathcal{I}} x_i} - \frac{D_i}{\sum_{i\in \mathcal{I}} D_i} \right| \leq K & $\forall i\in \mathcal{I}$ \label{eq:Const1_Orgut}\\
   S - \sum_{i\in \mathcal{I}} x_i - P = 0, & \label{eq:Const2_Orgut}\\
   x_i\leq C_i & $\forall i\in \mathcal{I}$ \label{eq:Const3_Orgut}\\
   x_i,P\geq 0 & $\forall i\in \mathcal{I}$, \label{eq:Const5_Orgut}
\end{subnumcases}
where $P$ represents the leftover food after the allocation of the donated supply $S$. We point out that $P$ represents Food-Bank Waste in Equation \eqref{eq:EFDef}. In this model, parameter $K$ is between 0 and 1 and controls the level of equity in the allocation of the total supply. Lower values of $K$ impose stricter (therefore higher) priority on equity and vice versa, wherein $K=0$ is the lower-bound of $K$ and corresponds to absolute equity. On the other hand, higher values of $K$ give more freedom to the objective function to choose solutions with higher efficiency. In order to find the upper-bound of $K$ for \cite{orgut2016modeling} model (denoted by $K^U$), we note that in perfect efficiency $P=0$. Therefore, from constraint \eqref{eq:Const2_Orgut}, we have $\sum_{i\in \mathcal{I}} x_i=S$. Thus, $K^U$ is the value of the objective function to the following mathematical model:

\begin{align}
Min \hspace{3pt}K^U  \label{eq:ObjFnc_Orgut_KU}
\end{align}\setcounter{equation}{11}
\begin{subnumcases}{s.t.}
   \left| \frac{x_i}{\sum_{i\in \mathcal{I}}x_i} - \frac{D_i}{\sum_{i\in \mathcal{I}} D_i} \right| \leq K^U & $\forall i\in \mathcal{I}$ \label{eq:Const1_Orgut_KU}\\
   \sum_{i\in \mathcal{I}}x_i  = S & \label{eq:Const2_Orgut_KU}\\  
   x_i\leq C_i & $\forall i\in \mathcal{I}$ \label{eq:Const3_Orgut_KU}\\
   x_i, K^U\geq 0 & $\forall i\in \mathcal{I}$. \label{eq:Const5_Orgut_KU}
\end{subnumcases}

Finally, to show the difference between the models more clearly, we divided the length between $\theta^L$ and $\theta^U$ as well as $0$ and $K^U$ to 50 pieces. Therefore, every point in Figures \ref{fig:EF-EQ_Gini} and \ref{fig:EQ-UEQ_Gini} corresponds to a value of $\theta$ ($K$) in our model (\cite{orgut2016modeling} model). We point out that some of the values of $\theta$ in our model yield identical equity and efficiency values, leading to reduced number of points for our model in the figures. Figures \ref{fig:Bin1-EF-EQ_Gini} and \ref{fig:Bin2-EF-EQ_Gini} show the efficient frontiers of our model and \cite{orgut2016modeling} using one minus the Gini coefficient to calculate the equity level. 

In both Figures \ref{fig:Bin1-EF-EQ_Gini} and \ref{fig:Bin2-EF-EQ_Gini}, the points resulting in perfect efficiency (Efficiency=1) correspond to $K=K^U$ in \cite{orgut2016modeling} model and $\theta=\theta^L$ in our model. On the other hand, the points resulting in perfect equity (Equity=1) correspond to $K=0$ in \cite{orgut2016modeling} model and $\theta=\theta^U$ in our model.

Our first observation from both Figures \ref{fig:Bin1-EF-EQ_Gini} and \ref{fig:Bin2-EF-EQ_Gini}is that while efficiency drops in an exponential manner as the value of $K$ decreases from $K^U$ in \cite{orgut2016modeling} model, the efficiency in our model drops in an almost linear manner while $\theta$ increases from its lower-bound $\theta^L$. This results in the fact that for the same levels of equity (efficiency), our model yields much higher levels of efficiency (equity). 

We note that when $K$ decreases from its upper-bound $K^U$, both equity and efficiency of \cite{orgut2016modeling} model decrease simultaneously. This is because constraint (\ref{eq:Const1_Orgut}) is not strictly capturing the concept of equity. Specifically, constraint \eqref{eq:Const1_Orgut} penalizes the deviations from $D_i/\sum_{i \in \mathcal{I}} D_i$, giving the priority of demand fulfillment to the agencies according to the magnitude of their demand. $D_i$ in constraint \eqref{eq:Const1_Orgut} can be interpreted as the weight (importance) of agency $i$. Therefore, agencies with larger demands will receive higher priority. When $K=K^U$, since constraint \eqref{eq:Const1_Orgut} is relaxed altogether, the agencies with smaller demands are treated equally as the ones with larger demands and receive some share of the supply, so equity is improved to some degree. However, once $K$ starts to decrease, constraint \eqref{eq:Const1_Orgut} forces the model to focus more on the agencies with larger demands, while the agencies with smaller demands are quickly neglected by the model. This leads to a decrease in equity. On the other hand, a decrease in $K$ inevitably leads to a decrease in efficiency by the definition of constraint \eqref{eq:Const1_Orgut}. Therefore, initially when $K$ decreases from its upper-bound $K^U$ both efficiency and equity decrease together.

As demand variability increases from Figure \ref{fig:Bin1-EF-EQ_Gini} to Figure \ref{fig:Bin2-EF-EQ_Gini}, the gap between the efficient frontiers of the two models increases. Specifically, while the variability increases, our model maintains almost the same equity-efficiency trade-off, making it robust against the variability in the data.

Regardless of the level of variations in the demand realizations, the two extremes in both models correspond to the cases with perfect equity (Equity = 1) and perfect efficiency (Efficiency = 1). At perfect equity both models result in the same efficiency. This is expected, because although the models have different objective functions, at perfect equity the solutions to the problems are unique, resulting in the same amount of waste at both models. However, at perfect efficiency the models may have infinitely many optimal solutions (please see our example in Section \ref{sec:ModAnalysis} corresponding to Definition \ref{defn: tightestlowbd}). In other words, since equity is not considered in the objective function under perfect efficiency, there may be infinitely many perfectly efficient solutions (Efficiency = 1), but with varying levels of equity. Among these solutions, our model chooses a dominant solution with considerable improvements in equity compared to that of \cite{orgut2016modeling}, especially under high demand variability.

\begin{figure}[H]
	\begin{center}
		\subfigure[\emph{low} demand variability]{\includegraphics[width=0.48\textwidth]{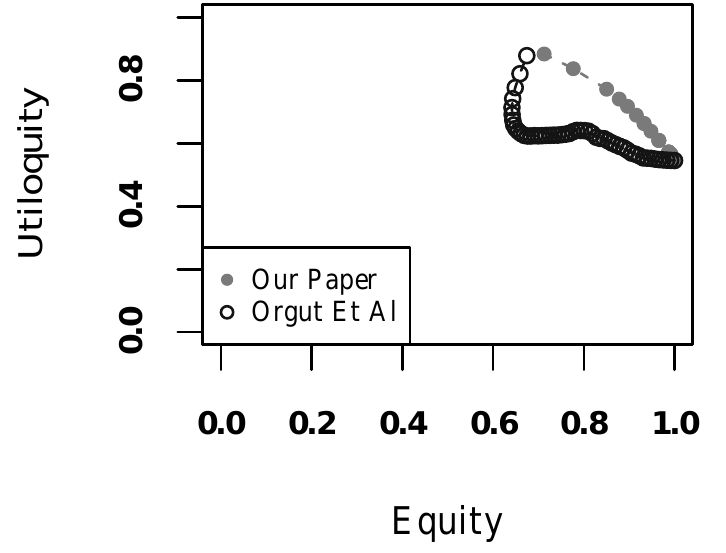} \label{fig:Bin1-EQ-UEQ_Gini}}
		\subfigure[\emph{high} demand variability]{\includegraphics[width=0.48\textwidth]{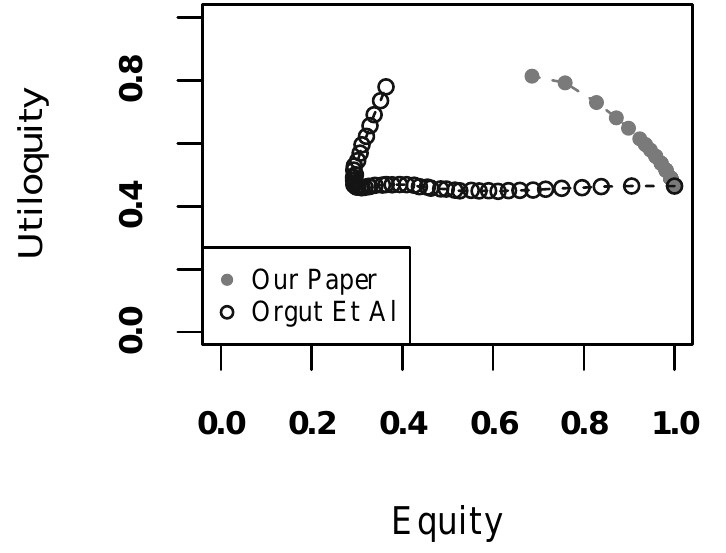} \label{fig:Bin2-EQ-UEQ_Gini}}
		\caption{Utilization Equity (Utiloquity) vs Equity using one minus Gini Coefficient as the equity measure.}
		\label{fig:EQ-UEQ_Gini}
	\end{center}
\end{figure}

In our setting, volunteer retention is directly proportional to utiloquity across the network. In other words, as the utiloquity increases, both over-utilization and under-utilization decrease across the agencies. Figure \ref{fig:EQ-UEQ_Gini} shows the trade-off between utiloquity and equity for both low and high demand variability. Similar to the case of equity-efficiency trade-off, we observed the same trend as the variability in demand increases in equity-utiloquity trade-off, i.e., the gain from our model compared to that of \cite{orgut2016modeling} increases. Additionally, although equity and utiloquity are conflicting objectives, for a given level of equity (utiloquity) our model achieves higher levels utiloquity (equity), regardless of variability level in the data. 

Our numerical experiments strongly suggest that, compared to our benchmark from the literature, our model not only achieves higher levels of efficiency for the same levels of equity, but also for a given efficiency it yields higher levels of equity both for the fill-rates and for the utilizations across the network. Therefore all three players in the food-bank supply chain simultaneously benefit: the food-bank, agencies, and the demand population.

To more concisely explain the equity-efficiency trade-off, in this section, we use a new measure called the price of equity. The price of equity in non-profit operations is defined as the required level of change in efficiency to achieve a unit of change in equity \citep{bertsimas2011price}. In other terms, price of equity is the difference of the efficient frontier curve (presented in Figure \ref{fig:EF-EQ_Gini}). 

Particular sources of data which are of interest to the managers of food-bank operations are the supply at the food-bank level and capacity at the agency level and the dynamics they create in the non-profit network operations. In this section, we quantify the effects of change in supply and capacity levels on the price of equity under two levels of variability in demand realizations and discuss some of our managerial insights. 

\begin{figure}[H]
	\centering
	\includegraphics[width=0.5\textwidth]{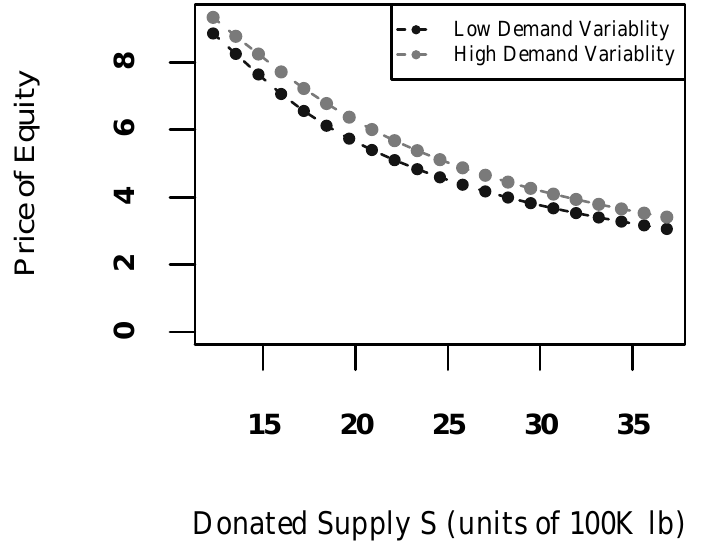}
	\caption{Price of Equity vs donated supply $S$ for both low and high demand variability.}
	\label{fig:PrEQ_S}
\end{figure}

Figure \ref{fig:PrEQ_S} shows the changes in price of equity as a function of available supply (denoted by $S$ in our model). As it is evident from Figure \ref{fig:PrEQ_S}, higher availability of supply reduces the amount of efficiency that should be sacrificed to achieve a single unit of increase in equity. In other words, the more a society becomes charitable, the easier it is for food-banks to become equitable in their distribution of food to the public. Another point of attention in Figure \ref{fig:PrEQ_S} is that the price of equity reduces with a diminishing rate as the supply increases, which indicates that initial level of increase in supply has the greatest impact in achieving a more equitable setting. 

Referring back to the definition of price of equity, we repeat that it shows the units of efficiency that need to be sacrificed in order to observe a unit change in equity. Figure \ref{fig:PrEQ_S} suggests that, at least at low supply, achieving equity is more difficult than achieving efficiency (price of equity is greater than one). In detail, price of efficiency can be defined as the reciprocal of the price of equity. Therefore a general observation from Figure \ref{fig:PrEQ_S} is that achieving higher equity requires more sacrifice from efficiency for the food-banks, especially when supplies are more scarce. 

Our final observation from Figure \ref{fig:PrEQ_S} is that as the variability in demand increases, the food-bank has to sacrifice more from its efficiency to achieve a single unit of equity for the same amount of supply. In other words, as the disparity of wealth in the society increases, it would be more difficult for food-banks to equitably distribute the supply they receive from the public.

\begin{figure}[H]
	\centering
	\includegraphics[width=0.5\textwidth]{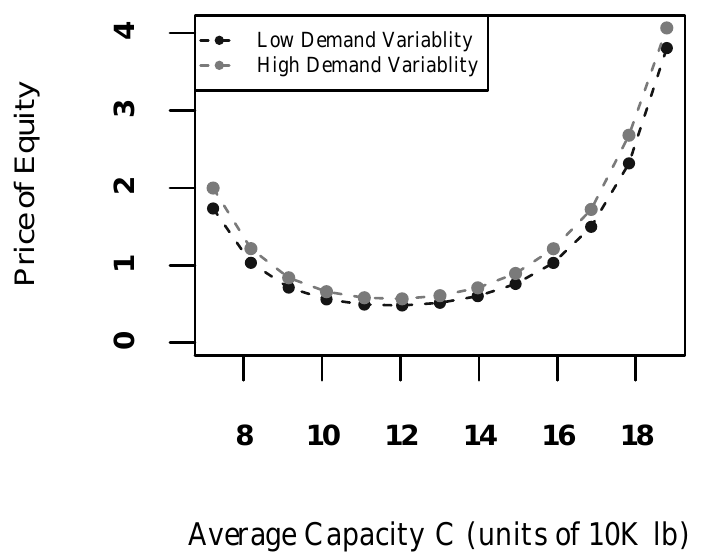}
	\caption{Price of Equity vs average capacity $C$ (across the agencies) for both low and high demand variability.}
	\label{fig:PrEQ_C}
\end{figure}

In traditional econometric models for capacity planning, the cost-capacity function has a unique minimizer. For instance, considering a traditional model for a warehousing system, cost of low capacity is the opportunity cost for the customer demands that could have been served, but are lost because there is not enough capacity for holding inventory. On the other hand, cost of high capacity is the idle space that is not used. 

Similarly, Figure \ref{fig:PrEQ_C} shows the changes in price of equity when average capacity increases in the network. Figure \ref{fig:PrEQ_C} suggests that there is a unique capacity at which price of equity is at its minimum. Specifically, a low capacity network of agencies in our setting behaves similar to having low supply which, as already discussed for Figure \ref{fig:PrEQ_S}, causes the price of equity to rise. This is intuitively correct as having low capacity to handle the donated food at the agencies limits the amount of supply they can receive and process. 

The cost of high capacity follows a similar behaviour. In detail, in our data from Feeding America, capacity is always a percentage of the demand in the network. This fact is a common phenomenon in practice as the capacity of the agencies typically falls short of their demand. Therefore, the cost of high capacity is due to the fact that an increase in capacities results in higher effective demands in the network (from Definition \ref{dfn: effective demand}). When effective demand increases and supply stays constant, it is as if the supply has become more scarce. Therefore, according to Figure \ref{fig:PrEQ_S}, the price of equity starts to rise after it meets its minimum in Figure \ref{fig:PrEQ_C}. 

In practice, Figure \ref{fig:PrEQ_C} demonstrates that contrary to the common belief, it is not always the case that expanding the capacity of food-banks will lead to ease in meeting higher equity levels. In other words, encouragement for volunteer involvement in food-bank operations and therefore enhancing their capacities, comes at a secondary level of priority when compared to reduction of poverty in the society. This is because, reducing poverty will lead to lower effective demands in the network and extends the decreasing behaviour of price of equity as the result of increase in capacities (before meeting the minimum) in Figure \ref{fig:PrEQ_C}. 

Our final observation from Figure \ref{fig:PrEQ_C} is that high variability in demand causes a shift in the price of equity as a function of capacity. In other words, expansion of capacities in food-banks leads to higher gains if the wealth is distributed more evenly in the society.

\section{Conclusions}
\label{sec:Conclusions}

Food-banks strive to distribute donations they receive from the public, government agencies, and grocery stores to the charitable agencies within their network in an equitable and efficient manner. In this paper, we present a new model to the food-bank donation allocation problem under equity and efficiency performance measures. Our model explicitly accounts for both efficiency and equity in the objective function and is capable of offering closed-form optimal solutions in perfect efficiency and perfect equity, while providing high quality solutions for the spectrum in between. 

Feeding America is one of the largest non-profit organizations, feeding the poverty population through a network of 200 food-banks across North America. Using the real data from Feeding America for their particular food-bank in North Carolina, we compare the results of our model against a benchmark from the literature in terms of equity, efficiency, and utilization equity (called by utiloquity in our paper). Our numerical study demonstrates considerable improvements in terms of all three performance measures simultaneously when compared to the benchmark. This means concurrently benefiting all three players in the food-bank supply chain: the food-bank, agencies, and the demand population.

Our sensitivity analysis demonstrates several interesting observations regarding the changes in price of equity when donated supply and agencies' capacities change. In particular, our sensitivity analysis demonstrates that the government should put its priority on helping the society to reduce poverty before investing on capacity expansions in charity organizations like food-banks. This will ensure that adding more capacity to the agencies associated with a food-bank will always lead to decreasing the price of equity. Additionally, we observed that encouraging the society towards charitability is always beneficial for the food-banks, however with a diminishing rate. Finally, our experiments demonstrate that decreasing disparity of wealth in the society will lead to lower costs for the food-banks to distribute public donations equitably to the poverty population.  

Our research can be extended in multiple directions. In this work, we considered a deterministic setting for our parameters of capacity and demand at the agencies. Consideration of uncertainty for either one of these parameters or both can perhaps lead the researchers to further realistic insights about the food-bank donation allocation problem. In another extension, we would like to refer the readers to the role of perishability in the decisions a food-bank makes. Particularly, food-banks typically receive close-to-expiry products from their donors and therefore integration of time into our model will further advance the insights in the science of decision-making in non-profit operations. 


\newpage
\begin{appendix}

\section{Proof of Theorem \ref{closed form of betaEQ}}\label{app: closedform betaEQ}

\noindent \underline{\emph{Proof}}

Here, we prove Theorem \ref{closed form of betaEQ}, the closed form of $\beta^{EQ}$.

First, if $\beta^a \leq \beta^b$, let $\boldsymbol{\beta}^{EQ} = (\beta^a, \cdots, \beta^a)$ and $\beta=\beta^a$. It is trivial to verify that ($\boldsymbol{\beta}^{EQ}$, $\beta$) is a feasible solution because $\beta^a \leq \beta^b \leq \tilde{\beta}_i$ for any $i \in \mathcal{I}$. On the other hand, there does not exists another $\beta^{EQ}>\beta^a$ but satisfies the perfect equity criterion, because $\beta^{EQ}D > \beta^aD = S$ is impossible. Hence, $\beta^{EQ}=\beta^a$.

Second, if $\beta^a > \beta^b$, then $\beta^{EQ}$ cannot be greater than $\beta^b$. If $\beta^{EQ}>\beta^b$ is true, then for any $j$ in the set $\mathcal{J}$, we have $\beta_j = \beta^{EQ} > \beta^b = \tilde{\beta}_j$ which violates the constraint (\ref{eq:Const2-efct}). Letting $\boldsymbol{\beta}^{EQ} = (\beta^b, \cdots, \beta^b)$ and $\beta=\beta^b$, it is trivial this solution is feasible. Since $\beta^b$ is the maximum possible $\beta^{EQ}$, $\beta^{EQ}=\beta^b$.

The above two claims complete the proof.

\begin{flushright}   Q.E.D.   \end{flushright}

\section{Proof of Theorem \ref{dominant efficient optimal}}\label{app: dominant efficient optimal}

\noindent \underline{\emph{Proof}}

Denote $\mathcal{I}_{i^*j^*}$ the first subset not fully filled, such that $\beta^{EF}_k D_k < \tilde{D}_k$ for any $k \in \mathcal{I}_{i^*j^*}$; other wise, if $S \geq \sum_{i \in \mathcal{I}} \tilde{D}_i$, then we have an unique optimal solution $\beta^{EF}_i = \tilde{D}_i/D_i$ ($\forall i \in \mathcal{I}$) which means all the agencies are filled at their maximum level -- the effective demand.

First, letting $\mathcal{B}$ be the set of all optimal solutions when $\theta=0$, we show that for any optimal solution $\bb=(\beta_1, \cdots, \beta_n)  \in \mathcal{B}$, $\beta_k=\beta^{EF}_k=0$ if $k \in \mI_{i}$ and $i>i^*$; $\beta_k=\beta^{EF}_k=\tilde{D}_i/D_i$ if $k \in \mI_{i}$ and $i<i^*$. In other words, the only difference part between the dominant solution $\bb^{EF}$ and other $\bb$'s in $\mathcal{B}$, if exists, is the fill-rates of agencies in the group $\mI_{i^*}$. Denoting $\mathcal{J}_1=\bigcup_{j=1}^{i^*-1} \mI_j$ and $\mathcal{J}_2=\bigcup_{j=i^*+1}^{n} \mI_j$, because $\beta^{EF}_k = \tilde{D}_i/D_i \geq \beta_k$ for any $k \in \mathcal{J}_1$ and $\beta^{EF}_k = 0 \leq \beta_k$ for any $k \in \mathcal{J}_2$, we have
\begin{align*}
 &\sum_{k\in\mathcal{J}_1}(\beta^{EF}_k-\beta_k)D_k + \sum_{k\in\mathcal{I}_{i^*}}(\beta^{EF}_k-\beta_k)D_k + \sum_{k\in\mathcal{J}_2}(\beta^{EF}_k-\beta_k) D_k = 0\\
\Rightarrow &  \sum_{k\in\mathcal{J}_1}(\beta^{EF}_k-\beta_k)D_k + \sum_{k\in\mathcal{I}_{i^*}}(\beta^{EF}_k-\beta_k)D_k = \sum_{k\in\mathcal{J}_2}(\beta_k-\beta^{EF}_k)D_k \\
\Rightarrow &  \sum_{k\in\mathcal{J}_1}(\beta^{EF}_k-\beta_k)D_k + D_{i^*}\sum_{k\in\mathcal{I}_{i^*}}(\beta^{EF}_k-\beta_k) = \sum_{k\in\mathcal{J}_2}(\beta_k-\beta^{EF}_k)D_k 
\end{align*}
The first equation comes from the observation $S < \sum_{i \in \mathcal{I}} \tilde{D}_i$. When there is no penalty of inequity, all the food supply will be shipped since the total effective demand is larger than the supply. On the other hand, notice that $\bb^{EF}, \bb \in \mathcal{B}$, so $\sum_{k\in \mI} (\beta^{EF}_k-\beta_k)=0$. However $\mI = \mathcal{J}_1 \cup \mI_{i^*} \cup \mathcal{J}_2$, we have $\sum_{k\in \mathcal{J}_1 \cup \mI_{i^*}} (\beta^{EF}_k-\beta_k)=\sum_{k\in \mathcal{J}_2} (\beta_k-\beta^{EF}_k)$. Therefore, if $\sum_{k\in \mathcal{J}_2} (\beta_k-\beta^{EF}_k) > 0$, then 
\begin{align*}
    \sum_{k\in \mathcal{J}_2} (\beta_k-\beta^{EF}_k) D_k \geq D_{i^*+1}\sum_{k\in \mathcal{J}_2} (\beta_k-\beta^{EF}_k) 
    & > D_{i^*} \sum_{k\in\mathcal{J}_1 \cup \mI_{i^*}} (\beta^{EF}_k-\beta_k) \\
    & > \sum_{k\in\mathcal{J}_1}(\beta^{EF}_k-\beta_k)D_k + D_{i^*}\sum_{k\in\mathcal{I}_{i^*}}(\beta^{EF}_k-\beta_k)
\end{align*}
gets contradiction. Hence, $\beta_k=\beta^{EF}_k=0$ if $k\in \mathcal{J}_2$. Furthermore, if $\sum_{k\in \mathcal{J}_1} (\beta^{EF}_k-\beta_k) > 0$, we have $\sum_{k\in \mathcal{J}_1} (\beta^{EF}_k-\beta_k)=\sum_{k\in \mathcal{I}_{i^*}} (\beta_k-\beta^{EF}_k)$ and
\begin{align*}
  D_{i^*}\sum_{k\in\mathcal{J}_1}(\beta^{EF}_k-\beta_k) > \sum_{k\in\mathcal{J}_1}(\beta^{EF}_k-\beta_k)D_k = D_{i^*}\sum_{k\in\mathcal{I}_{i^*}}(\beta_k-\beta^{EF}_k)
\end{align*}
gets contradiction. Hence, $\beta_k=\beta^{EF}_k=\tilde{D}_i/D_i$ if $k\in \mathcal{J}_1$ and $\sum_{k\in \mathcal{I}_{i^*}} (\beta_k-\beta^{EF}_k)=0$.

Now, based on the results above, we have $\beta_k = \beta^{EF}_k$ when $k\in \mI/\mI_{i^*}$. On the other hand, by how Algorithm \ref{alg:domiEfficient} (Greedy method) constructs $\bb^{EF}$, it is obvious $\max_{k\in \mI_{i^*}} \beta_k \geq \max_{k\in \mI_{i^*}} \beta^{EF}_k$. Therefore, $\beta \geq \beta^{EF}$. When $\theta>0$, consider the difference between objective function values as follows.
\begin{align*}
    & \sum_{k \in \mathcal{I}} [\beta^{EF}_k-\theta(\beta^{EF}-\beta^{EF}_k)] - \sum_{k \in \mathcal{I}} [\beta_k-\theta(\beta-\beta_k)]\\
  =&\theta |\mI|(\beta-\beta^{EF}) + (1+\theta) \sum_{k \in \mI} (\beta^{EF}_k-\beta_k) \\
  =&\theta |\mI|(\beta-\beta^{EF}) \geq 0
\end{align*}
completes the proof.

\begin{flushright}   Q.E.D.   \end{flushright}

\section{Performance of the Models Under Other Measures of Inequity}\label{app: EQ_OtherMeasures}
In this section, we report the comparison between the two models based on the remaining four measures of inequity: 

    \begin{align}
        &\text{Coefficient of Variation} =\frac{\sqrt{\frac{\sum_{i \in \mathcal{I}}(E_i - \bar{E})^2}{n}}}{\bar{E}}  \\  
        &\text{Variance} =\frac{\sum_{i \in \mathcal{I}}(E_i - \bar{E})^2}{n}  \\ 
        &\text{Mean Absolute Deviation} =\frac{\sum_{i \in \mathcal{I}}\left|E_i - \bar{E}\right|}{n}  \\ 
        &\text{Range} = E_{max} - E_{min}
    \end{align}

In order to draw the efficient frontier we need to convert the measures to return equity. However, while Gini Coefficient returns a value between 0 and 1, the rest of the above-mentioned four measures return a non-negative value which is not necessarily between 0 and 1. Therefore, similar to \cite{orgut2016modeling}, we first standardize our data for each point in each measure $\boldsymbol{E}=(E_1,E_2, \dots, E_n)$ to be between 0 and 1 in its standard form $\boldsymbol{\hat{E}}=(\hat{E}_1,\hat{E}_2, \dots, \hat{E}_n)$ according to the following:
    \begin{equation}
        \hat{E}_i=\frac{E_i - E_{Min}}{E_{Max} - E_{Min}}
        \end{equation}

Our equity vector then is calculated as $1-\boldsymbol{\hat{E}}$.

\begin{figure}[H]
	\begin{center}
		\subfigure[\emph{low} demand variability]{\includegraphics[width=0.48\textwidth]{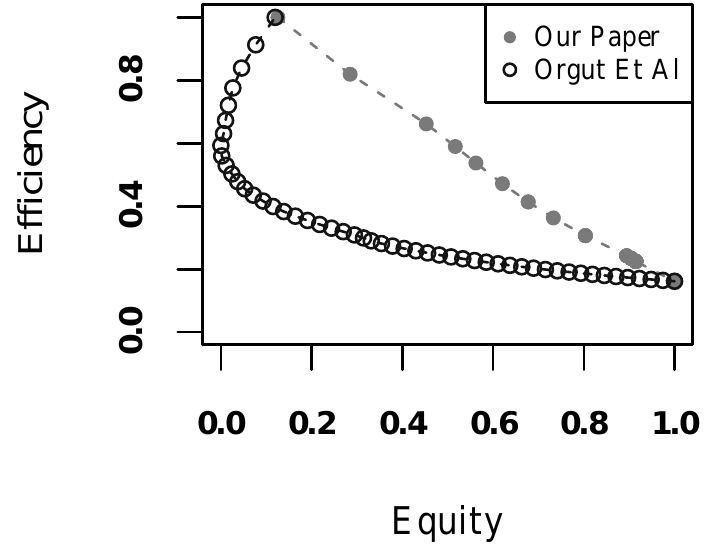} \label{fig:Bin1-EF-EQ_CV}}
		\subfigure[\emph{high} demand variability]{\includegraphics[width=0.48\textwidth]{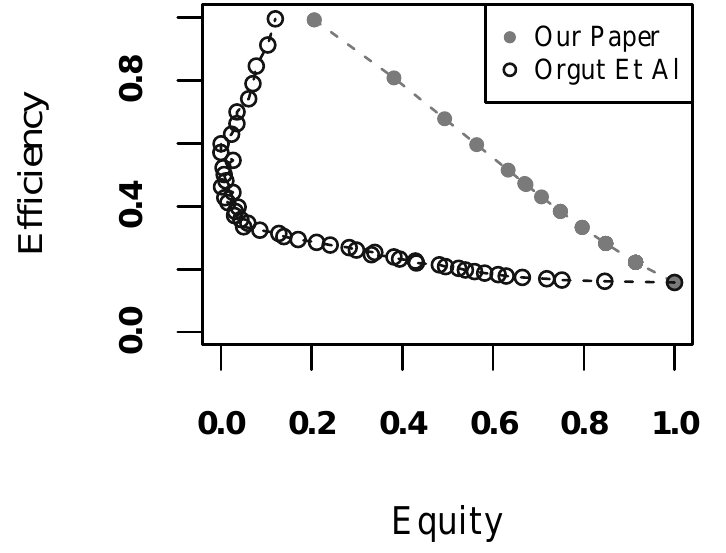} \label{fig:Bin2-EF-EQ_CV}}
		\caption{Efficient frontier using \emph{Coefficient of Variation} as the inequity measure.}
		\label{fig:EF-EQ_CV}
	\end{center}
\end{figure}

\begin{figure}[H]
	\begin{center}
		\subfigure[\emph{low} demand variability]{\includegraphics[width=0.48\textwidth]{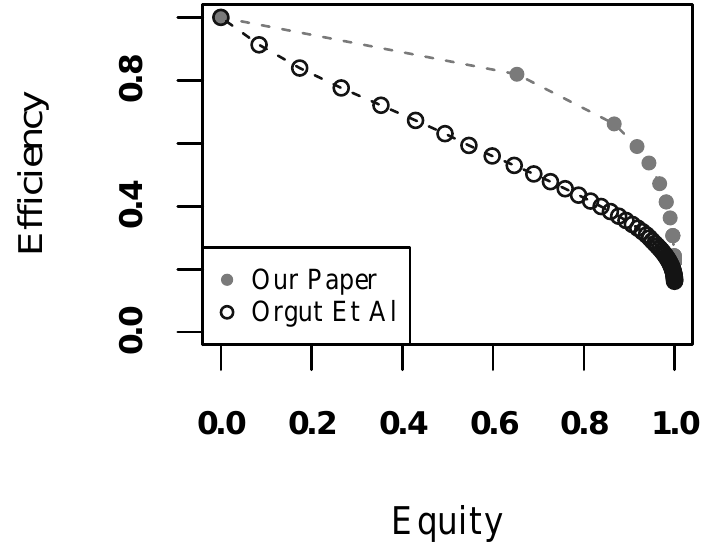} \label{fig:Bin1-EF-EQ_Var}}
		\subfigure[\emph{high} demand variability]{\includegraphics[width=0.48\textwidth]{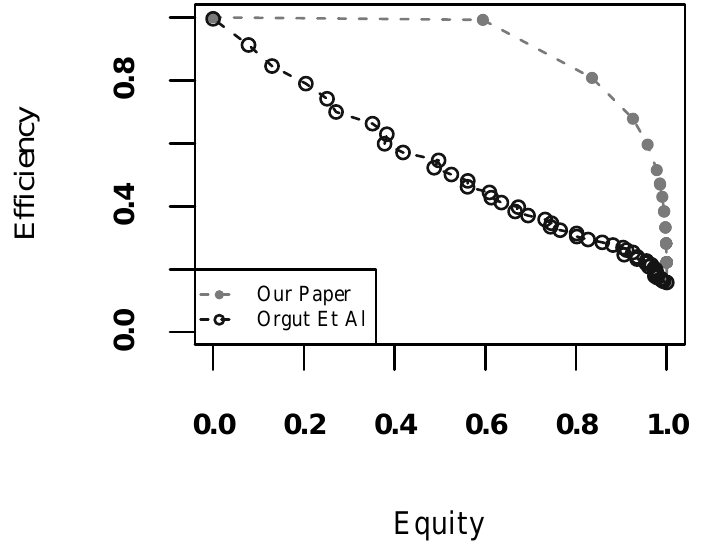} \label{fig:Bin2-EF-EQ_Var}}
		\caption{Efficient frontier using \emph{Variance} as the inequity measure.}
		\label{fig:EF-EQ_Var}
	\end{center}
\end{figure}

\begin{figure}[H]
	\begin{center}
		\subfigure[\emph{low} demand variability]{\includegraphics[width=0.48\textwidth]{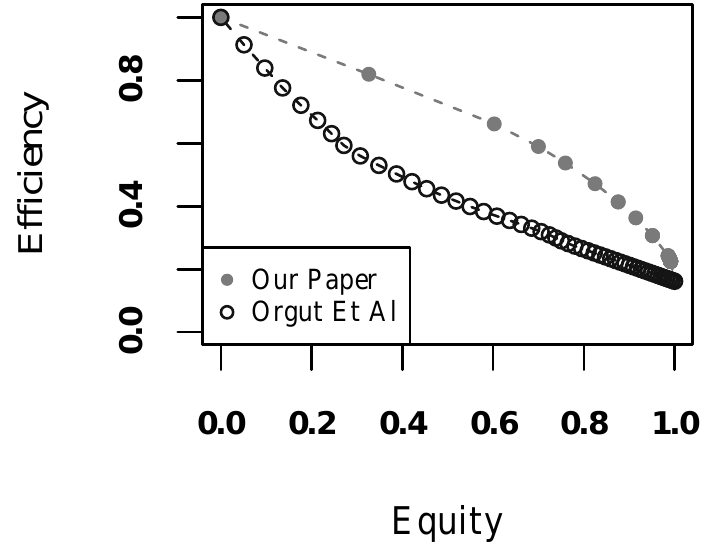} \label{fig:Bin1-EF-EQ_MAD}}
		\subfigure[\emph{high} demand variability]{\includegraphics[width=0.48\textwidth]{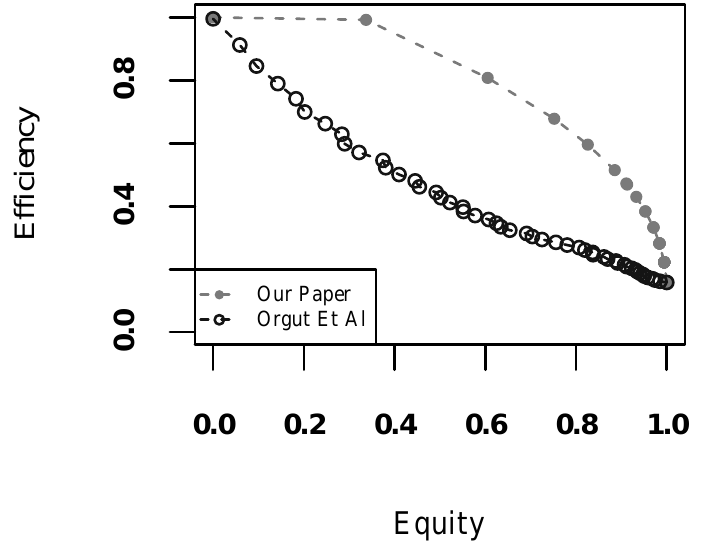} \label{fig:Bin2-EF-EQ_MAD}}
		\caption{Efficient frontier using \emph{Mean Absolute Deviations} as the inequity measure.}
		\label{fig:EF-EQ_MAD}
	\end{center}
\end{figure}

\begin{figure}[H]
	\begin{center}
		\subfigure[\emph{low} demand variability]{\includegraphics[width=0.48\textwidth]{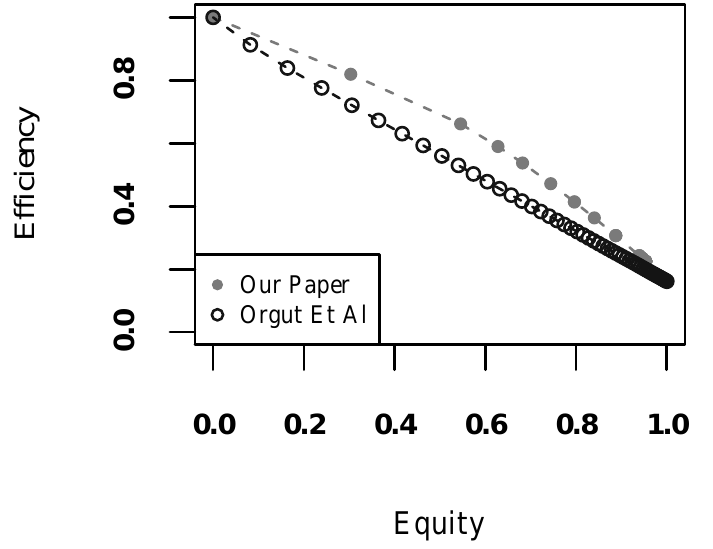} \label{fig:Bin1-EF-EQ_Range}}
		\subfigure[\emph{high} demand variability]{\includegraphics[width=0.48\textwidth]{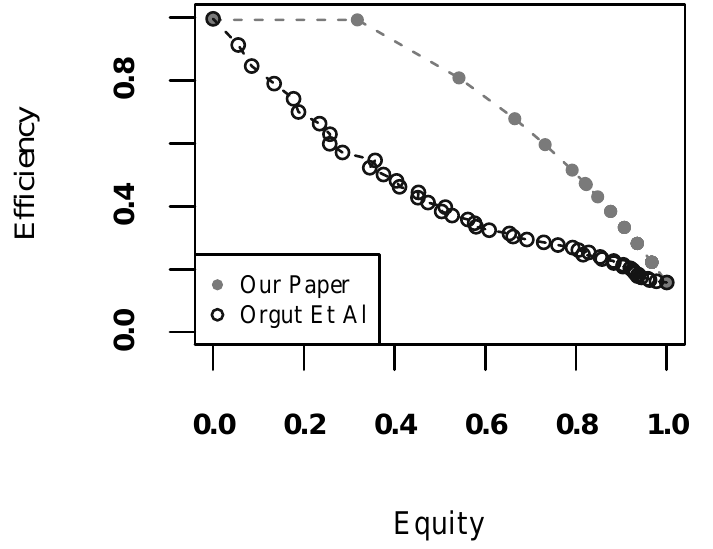} \label{fig:Bin2-EF-EQ_Range}}
		\caption{Efficient frontier using \emph{Range} as the inequity measure.}
		\label{fig:EF-EQ_Range}
	\end{center}
\end{figure}

\end{appendix}

\newpage
\bibliographystyle{ormsv080}

\bibliography{main}

\end{document}